\documentclass{article}
\usepackage{psfrag}
\usepackage{graphicx}
\usepackage{subfigure}
\usepackage{amsmath,amssymb,amsfonts}

\newtheorem{thrm}{Theorem}[section]
\newtheorem{prpstn}[thrm]{Proposition}
\newtheorem{lmm}[thrm]{Lemma}
\newtheorem{dfntn}[thrm]{Definition}

\newtheorem{rmrk}[thrm]{Remark}

\newtheorem{crllr}[thrm]{Corollary}

\def\sqw{\hbox{\rlap{\leavevmode\raise.3ex\hbox{$\sqcap$}}$%
\sqcup$}}
\def\findem{\ifmmode\sqw\else{\ifhmode\unskip\fi\nobreak\hfil
\penalty50\hskip1em\null\nobreak\hfil\sqw
\parfillskip=0pt\finalhyphendemerits=0\endgraf}\fi}

\def\dem {\noindent {\bf Proof: }}

\def \NNN {\mathcal{N}}
\def \CCC {\mathcal{C}}
\def \itt {{n}}
\def \CC { A} 
\def \BB { B} 
\def \Dvec {D}
\def \ee {\mathrm{e}}
\def \hh {{\bf h}}

\def \qqq {{\bf q}}
\def \qq {\mathrm{q}}
\def \uu {{\bf u}}
\def \nn {{\bf n}}
\def \UU {{\bf U}}
\def \WW {{\bf W}}

\def \GG {{\bf G}}
\def \ww {{\bf w}}
\def \vv {{\bf v}}
\def \xx {{\bf x}}
\def \yy {{\bf y}}
\def \zz {{\bf z}}
\def \Gmat {G}

\def \R {\mathbb{R}}
\def \NN {\mathrm{N}}
\def \Qc {\tilde{Q}}
\def \P {\textsc{P}}

\def \virg {\, , \,\,}
\def \dsp {\displaystyle}

\def \bflambda {{\boldsymbol{\lambda}}}
\def \bfmu {\boldsymbol{\mu}}

\def \vsp {\vspace{6pt}}        

\date{ }
\title{A discrete contact model for crowd motion}
\author{Bertrand Maury\\ Universit\'e de Paris-Sud\\UMR du CNRS 8628\\F-91405 Orsay
Cedex\\bertrand.maury@math.u-psud.fr \and Juliette Venel\\ Universit\'e de Paris-Sud\\UMR du CNRS 8628\\F-91405 Orsay
Cedex\\juliette.venel@math.u-psud.fr}

\begin{document}

%
\maketitle

%
%
\subsection*{Abstract}
The aim of this paper is to develop  a crowd motion model designed to handle highly packed situations. 
The  model we propose  rests on two principles:
We first define a
spontaneous velocity which corresponds to the velocity each
individual would like to have in the absence of other people; The
actual velocity is then computed as the projection of the spontaneous 
velocity onto the set of admissible velocities (i.e. velocities which
do not violate the non-overlapping constraint). 
We describe here the underlying mathematical framework, and we explain
how recent results by J.F. Edmond and L. Thibault on the sweeping process by uniformly prox-regular sets can be adapted to handle this situation in terms of well-posedness. 
We propose a numerical scheme for this contact dynamics model, based on a prediction-correction algorithm.
Numerical illustrations are finally presented and discussed.

\subsection*{R\'esum\'e}
Nous proposons un mod\`ele de mouvements de foule orient\'e vers la gestion de  configurations tr\`es denses. Ce mod\`ele repose sur deux principes: tout d'abord nous d\'efinissons une vitesse souhait\'ee correspondant \`a la vitesse que les individus aimeraient avoir en l'absence des autres; la vitesse r\'eelle est alors obtenue comme projection de la vitesse souhait\'ee sur un ensemble de vitesses admissibles (i.e. qui respectent la contrainte de non-chevauchement). Nous d\'ecrivons le cadre math\'ematique sous-jacent et nous expliquons comment certains r\'esultats de J.F. Edmond et L. Thibault sur les processus de rafle par des ensembles uniform\'ement prox-r\'eguliers peuvent \^etre utilis\'es pour prouver le caract\`ere bien pos\'e de notre mod\`ele.
Nous proposons un sch\'ema num\'erique pour ce mod\`ele de dynamique des contacts bas\'e sur un algorithme de type pr\'ediction-correction. Enfin des r\'esultats num\'eriques sont pr\'esent\'es et comment\'es.

%
%
%
\section*{Introduction}

Walking behaviour of pedestrians has given rise to a large amount of empirical studies over the last decades.
Qualitative data (preferences, walk tendencies) have been collected by Fruin~\cite{Fruin}, Navin,  Wheeler~\cite{Navin}, Henderson~\cite{Henderson} and, more recently, by Weidmann~\cite{Weidmann}. 
From these observations, several strategies for crowd motion modelling have been proposed, 
and can be classified with respect  to the way they handle people density (Lagrangian description of individuals or macroscopic approach), and to the nature of motion phenomena (deterministic or stochastic).
Among discrete and stochastic models, let us mention 
Cellular Automata~\cite{Blue,Burstedde,Nagel,Schad}, models based on networks~\cite{Exodus} as route choice models~\cite{BorgersTim0, BorgersTim} and queuing models~\cite{Lovas,Yuhaski}.
In these models, each cell or node is either empty or occupied by a single person and people's motion always satisfies this rule. In cellular automata models, there are two manners of moving people during a time step. With the first one, positions are updated one by one with a random order (\textit{Random Sequential Update}). The second method consists of updating simultaneously all positions (\textit{Parallel Update}). If several people want to reach the same cell, only one of them (randomly chosen) is allowed to move. 
In route choice models, people move on a network. Each model is based on a route choice set. Most choice set generation procedures are based on shortest route search and use shortest paths algorithms.
Queuing models use Markov-chain models to describe how pedestrians move from one node of the network to another. 

In~\cite{Helbsocforce}, a microscopic model called social force model is presented. It describes crowd motion with a system of differential equations. The acceleration of an individual is obtained according to Newton's law. Several forces are introduced as for example a term describing the acceleration towards the desired velocity or a repulsion force reflecting that a pedestrian tends to keep a certain distance from other people and obstacles.
 Moreover macroscopic models have been proposed. In~\cite{Henderson}, pedestrian traffic dynamics is firstly compared with fluid dynamics. Some models~\cite{Helbfluid,Henderson,Hooggas} are based on gas-kinetic theory. Other models~\cite{Hoog2,Hoog1,Hughes,Hughes2} rest on a set of partial differential equations describing the conservation of flow equation.  
 
 Several softwares have been developped: PedGo~\cite{PedGo}, SimPed~\cite{Daamenthese}, Legion~\cite{Still}, Mipsim~\cite{Hooggas} or Exodus~\cite{Exodus}.
Some commonly observed  collective patterns are now considered  as standard benchmarks for those numerical simulations. 
 %
Among  these phenomena of self-organization, there is the formation of lanes formed naturally by people moving in opposite directions.
In this way, strong interactions with
oncoming pedestrians are reduced, and a higher walking speed is possible. Another phenomena is the formation of arches upstream the  exit during the 
  evacuation of a room. These patterns are recovered by CA-models~\cite{Kirchevac, Schadant} and by the social force model~\cite{Helbsocforce, Helbpanic}. 

The case of evacuation in emergency situations is of particular importance in terms of applications  (observance of security rules, computer-assisted design of public buildings, appropriate positionning of   exit signs). 
Numerical simulations may allow to  estimate evacuation time (to be compared for example with the duration of fire propagation) and also to predict areas where high density will appear. 
As pointed out by Helbing~\cite{Helbpanic}, emergency situations do not fit into the standard framework of pedestrian traffic flow.
When people stroll around without hurry, they tend to keep a certain distance from each other and from obstacles. 
In an emergency situation, the motion of individuals is governed by different rules. In particular, the contact with walls or other people is no longer avoided. 
Some strategies have been proposed to adapt social walk models to highly congested situation (see again~\cite{Helbpanic}). 
We propose here an approach which relies  on the very consideration that actual motion in emergency situations is governed by the opposition between achievement of individual satisfaction (people struggle to escape as quickly as possible, regardless of the global efficiency) and congestion.  
In particular, we aim  at integrating the direct conflict between people in the model, in order to estimate in some way interaction forces between them, and therefore provide a way to estimate the local risk of casualties.

The microscopic model we propose rests on two principles. On the one hand, each individual has a spontaneous velocity that he would like to have in the absence of other people. On the other hand, the actual velocity must take into account congestion. Those two principles lead us to define the actual velocity field as the projection of the spontaneous velocity onto the set of admissible velocities (regarding the non-overlapping constraints). The flexibility of this model lies in its first point: every choice of spontaneous velocity can be made and so every existing model for predicting crowd motion can be integrated here. The key feature of the model is the second point which concerns  handling of contacts.

By specifying the link between these two velocities, the evolution problem takes the form of a first order differential inclusion. 
This type of evolution problem has been extensively  studied in the 1970's, 
 with the theory of maximal monotone operators  (see e.g.~\cite{Brezis}). A few years later, J.J. Moreau considered similar problems with time-dependent multivalued operator, namely sweeping processes by convex sets (see~\cite{Moreausweep}).
 Since then, important improvements have been developped by weakening the convexity assumption with the concept of prox-regularity. The well-posedness of our evolution problem can be established by means of recent results of J.F. Edmond and L. Thibault~\cite{Thibrelax} concerning sweeping processes by uniformly prox-regular sets. 

The paper is structured as follows: In Section 2, we present the model and establish its well-posedness; In Section 3, we propose a numerical scheme, and detail the overall solution method. 
Section 4 is devoted to some illustrations of the numerical algorithm.

\section{Modelling}
\begin{figure}
\centering
\psfrag{r}[l]{$r$}
\psfrag{e}[l]{$\ee_{ij}(\qqq)$}
\psfrag{me}[l]{$\hspace{-0.5cm}\vspace{-2cm} -\ee_{ij}(\qqq)$}
\psfrag{qj}[l]{$\qq_j$}
\psfrag{qi}[l]{$\qq_i$}
\psfrag{d}[l]{ $\hspace{-0.5cm} D_{ij}(\qqq)$}
\includegraphics[width=0.4\textwidth]{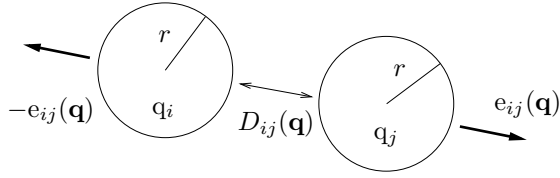}
\caption{Notations.}
\label{fig:notations}
\end{figure}

We consider $N$ persons identified to rigid disks. For convenience, the disks are supposed here to have the same radius $r$. The
centre of the i-th disk is denoted by $\qq_i$ (see
Fig.~\ref{fig:notations}). 
Since overlapping is forbidden, the  
vector of positions $ \dsp \qqq=(\qq_1,..,\qq_N) \in
\R^{2N}  $ (equipped with the euclidean norm) is required to  belong to the following set: 
\begin{dfntn}[Set of feasible configurations]
$$\dsp Q=\left\{ \qqq  \in
\R^{2N} \virg  D_{ij}( \qqq) \geq 0\quad   \forall \,i < j  \right\}, $$ 
where $ D_{ij}( \qqq)=|\qq_i-\qq_j
  |-2r $ is the signed distance between disks $i$ and $j$. 
\label{def:Q0}
\end{dfntn}
We consider as given the vector of spontaneous velocities denoted by $$ \dsp \UU(\qqq)=( \UU_1(\qqq),\ldots,  \UU_N(\qqq)) \in \R^{2N}. $$
$\UU_i$ is the spontaneous velocity of individual $i$, which may depend on its own position ($\UU_i= \UU_i(\qqq_i)$, see Section~\ref{sec:resnum} for examples of such a situation), but also on other people's positions, that is why we keep here $\UU_i= \UU_i(\qqq)$.
To define the actual velocity, we introduce the following set:
\begin{dfntn} {(Set of feasible velocities)} \mbox
\newline $$ \dsp \CCC_{\qqq}=\left\{  \vv  \in
\R^{2N} , \ \forall i<j \hspace{5mm} D_{ij}(\qqq)=0 \hspace{3mm}
\Rightarrow \hspace{3mm}
\GG_{ij}(\qqq)\cdot \vv \geq 0 \right\},  $$ with
$$\dsp \GG_{ij}(\qqq)=\nabla D_{ij}(\qqq)=
(0,\dots,0, -\ee _{ij}(\qqq),0,\dots,0,\ee _{ij}(\qqq),0,\dots,0 )\in
\R^{2N} \textmd{  and  } \ee _{ij}(\qqq)= \dsp \frac{\qq_j-\qq_i}{|\qq_j-\qq_i|}.
 $$
\label{def:Cq}
\end{dfntn}
 The actual velocity field is defined as 
 the feasible field which is the closest to $\UU$ 
 in the least square sense, which writes 
\begin{equation}
\left\{
\begin{array}{l}
 \dsp \frac{d\qqq}{dt} = \P_{\CCC_{\qqq}} \UU(\qqq),
\vspace{2mm} \\
\qqq(0) = \qqq_0 \in Q,
\end{array}
\right.
\label{eq:projCq}
\end{equation}
where $ \P_{\CCC_{\qqq}} $ denotes the euclidean projection onto the
closed convex cone $\CCC_{\qqq}$. 
\begin{rmrk}
Despite its formal simplicity, this model does not fit directly into a
standard framework. Indeed the set $\CCC_\qqq$ does not continuously depend on $\qqq$. If no contact holds, the velocity is not constrained and $\CCC_\qqq= \R^{2N} $. With a single contact, the set $\CCC_\qqq$ becomes a half-space.
\end{rmrk}
\section{Mathematical framework}
\subsection{Reformulation}
\begin{figure}
\centering
\psfrag{e}[c]{$\eta$}
\psfrag{x}[c]{$\xx$}
\psfrag{s}[c]{$S$}
\includegraphics[width=0.23\textwidth]{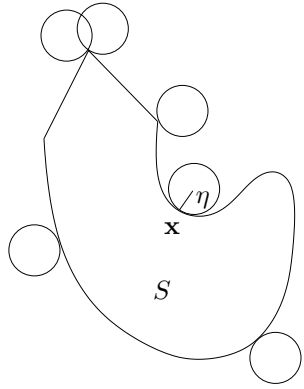}
\caption{$\eta $-prox-regular set.}
\label{fig:ensproxreg} 
\end{figure}

Let us reformulate the problem by introducing $\NNN_\qqq $, the outward
normal cone to the set of feasible configurations $Q$, which is
defined as the polar cone of $\CCC_\qqq $.
\begin{dfntn} {(Outward normal cone)} \mbox 
 \newline $$
\NNN_\qqq = \CCC_\qqq^\circ
=  \left\{\ww \in \R^{2N}\virg \ww \cdot \vv \leq 0 \quad \forall \vv \in \CCC_\qqq \right\}.
$$
\label{def:Nq}
\end{dfntn}
\begin{figure}
\centering
\psfrag{a}[l]{$D_{12}<0 $}
\psfrag{b}[l]{$D_{13}<0 $}
\psfrag{c}[l]{$D_{34}<0 $}
\psfrag{Q}[l]{$ Q_0 $}
\psfrag{u}[l]{$\qqq_2$}
\psfrag{v}[l]{$\qqq_1$}
\psfrag{e}[l]{$\NNN_{\qqq_2} $}
\psfrag{f}[l]{$\CCC_{\qqq_2} $}
\psfrag{d}[l]{$\NNN_{\qqq_1} $}
\psfrag{g}[l]{$\CCC_{\qqq_1} $}
\includegraphics[width=0.55\textwidth]{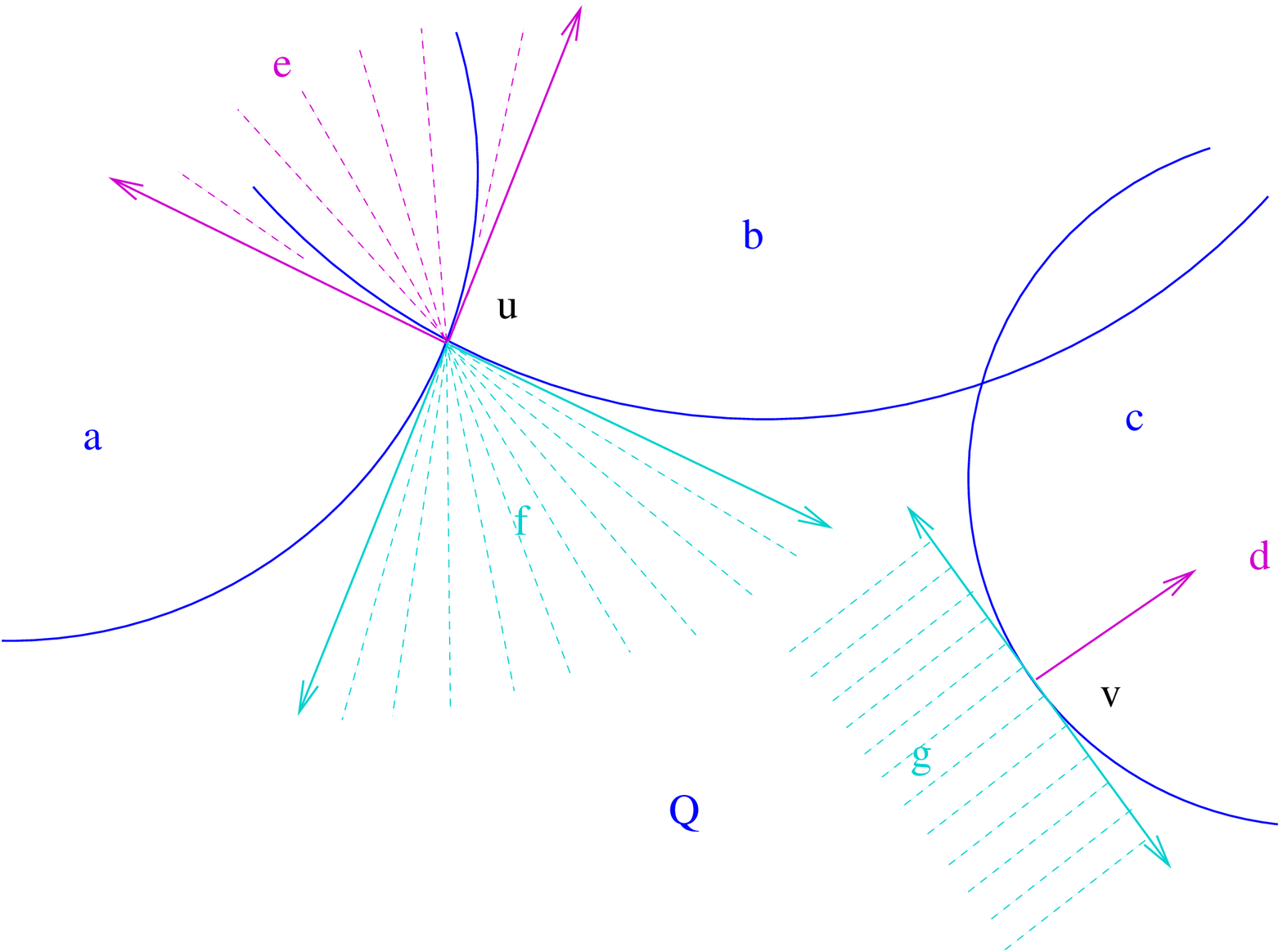}
\caption{Cones $ \CCC_\qqq$ and $\NNN_\qqq$.}
\label{fig:NqCq}
\end{figure}
\begin{rmrk}
In Figure~\ref{fig:NqCq}, we represent the set $Q \subset \R^{2N}$ which is defined as an intersection of convex sets' complements. In the case of a single contact (configuration $\qqq_1$), we remark that the cone $ \NNN_{\qqq_1}$ is generated by the vector $-\GG_{34}(\qqq_1)$ that is up to a constant, the outward normal vector to the domain $D_{34}\geq 0$. In the case of two or more contacts, the configuration $\qqq_2$ does not belong to a smooth surface and the cone $\NNN_{\qqq_2}$ (generated by $-\GG_{12}(\qqq_2) $ and $ -\GG_{13}(\qqq_2)$) generalizes somehow the notion of the outward normal direction.  
\label{denoNq}
\end{rmrk}
Thanks to Farkas' Lemma (see~\cite{Ciarlet}), the outward normal cone can be expressed
\begin{equation}
\NNN_\qqq = \left\{-\sum \lambda_{ij} \GG_{ij}(\qqq) \virg \lambda_{ij} \geq 0 \virg
D_{ij}(\qqq) > 0 \Longrightarrow \lambda_{ij} = 0
\right\} .
\label{exprNq}
\end{equation}
Let us recall the classical orthogonal decomposition of a Hilbert space as the sum of mutually polar cone (see~\cite{Moreaucones}) :
\begin{equation}
\P_{\CCC_\qqq} + \P_{\NNN_\qqq} = Id.
\label{decomp}
\end{equation}
Using this property, we get:
\begin{equation}
\frac {d \qqq}{dt} = \P_{\CCC_\qqq} \UU (\qqq)= \UU(\qqq) - \P_{\NNN_\qqq} \UU(\qqq).
\label{eqdiff}
\end{equation}
Since $\P_{\NNN_\qqq} \UU(\qqq) \in \NNN_\qqq$, we obtain a new formulation for (\ref{eq:projCq}) 
\begin{equation}
\left\{
\begin{array}{l}
 \dsp \frac{d\qqq}{dt} + \NNN_\qqq \ni \UU(\qqq),
\vspace{2mm} \\
\qqq(0)=\qqq_0.
\end{array}
\right.
\label{incldiff}
\end{equation}
The problem reads as a first order differential inclusion involving the multivalued operator $\NNN$.
\begin{rmrk}
In the absence of contacts in the configuration $\qqq$, the set of feasible velocities $\CCC_\qqq$ is equal to the whole space $\R^{2N}$, and consequently the outward normal cone $\NNN_\qqq$ is reduced to $\{0\}$. In that case, the first relation of $($\ref{incldiff}$)$ states that the actual velocity equals to the spontaneous velocity: $$\dsp \frac{d\qqq}{dt} = \UU(\qqq) .$$ If any contact exists, the differential inclusion means that the configuration $\qqq$, submitted to $\UU(\qqq)$, has to evolve while remaining in $Q$. 
\label{interpretationincludif}
\end{rmrk}
%
 Let us first study a special situation where standard theory can be applied. 
Consider $N$ individuals in a corridor.
In that case, as people cannot leap accross each other, it is natural to restrict the set of feasible configurations to one of its connected components:
$$
Q = \{ \qqq=(\qq_1,\dots,\qq_N) \in \R^N\virg \qq_{i+1}-\qq_i \geq 2r
\} .
$$
In this very situation, as $Q$ is closed and convex,  the
multivalued operator $\qqq \longmapsto \NNN_\qqq$ identifies to the subdifferential of the indicatrix function of $Q$:
$$
  \partial I _{Q}(\qqq) =  \{ \vv ,   I _{Q}(\qqq) + (\vv,\hh) \leq  I _{Q}(\qqq + \hh) \quad \forall \hh
\}
\virg
 I _{Q}(\qqq)= \left |
\begin{array}{rcl} 0&\hbox{ if }& \qqq \in  {Q} \\
+\infty&\hbox{ if }& \qqq \notin  {Q}
\end{array}\right .
$$
therefore $\qqq\longmapsto \NNN_\qqq$ is maximal monotone.
In that case, as soon as the spontaneous velocity is regular (say Lipschitz), 
standard theory (see e.g. Brezis~\cite{Brezis}) ensures well-posedness.
%
\begin{figure}
\centering
\psfrag{qqq2}[c]{$\bar{\qq}_2$}
\psfrag{qqq1}[c]{$\bar{\qq}_1$}
\psfrag{qq2}[c]{$\tilde{\qq}_2$}
\psfrag{qq1}[c]{$\tilde{\qq}_1$}
\psfrag{q2}[c]{$\qq_2$}
\psfrag{q1}[c]{$\qq_1$}
\psfrag{c1}[c]{\small{configuration} $\qqq$}
\psfrag{c2}[c]{\small{configuration} $\widetilde{\qqq}$}
\psfrag{c3}[c]{\small{configuration} $ \dsp \bar{\qqq}=\frac{\qqq+\widetilde{\qqq}}{2}$}
\includegraphics[width=0.75\textwidth]{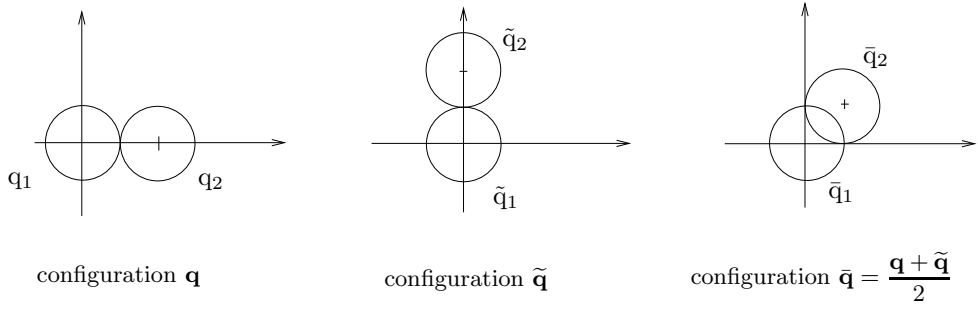}
\caption{Lack of convexity.}
\label{fig:nonconvex}
\end{figure}
Yet, as illustrated in Figure~\ref{fig:nonconvex}, 
$Q$ is not convex in general and the operator $\qqq\longmapsto
\NNN_\qqq$ is not monotone. So we cannot apply the same arguments as in the case of a straight motion. By lack of convexity, the projection onto $Q$ is not everywhere well-defined.
 However the set $Q$ satisfies a weaker property in the sense that the projection onto $ Q$ is still well-defined in its neighbourhood. Indeed,  $Q$ is uniformly prox-regular, which is the suitable property to
ensure well-posedness. Let us give some definitions to
specify the general mathematical framework.
\begin{dfntn}
Let $S$ be a closed subset of a Hilbert space $H$.\\ We define the proximal
normal cone to $S$ at $\xx$ by:  $$ \NN(S,\xx)=
\left \{
\vv \in H, \  \exists \alpha >0, \  \xx \in \P_S(\xx + \alpha
\vv) \right \}, $$
where  $$\P_S(\yy)=\{\zz \in S, \ d_S(\yy)= |\yy -\zz |\}, \textmd{ with } d_S(\yy)=\underset{\zz \in S}{\mathrm{min} }|\yy -\zz |. $$
\end{dfntn}
Following~\cite{Clarke}, we define the concept of uniform prox-regularity as follows:
\begin{dfntn}
Let $S$ be a closed subset of a Hilbert space $H$. $S$ is said $\eta$-prox-regular if for
all $\xx \in \partial S $ and $\vv \in \NN(S,\xx), \  |\vv|=1 $ we have:
$$ B(\xx+\eta \vv, \eta)\cap S = \emptyset.  $$
\label{caracpr}
\end{dfntn}
In an euclidean space, $S$ is $\eta$-prox-regular if an external tangent ball with radius smaller than $\eta$ can be rolled around it (see Fig~\ref{fig:ensproxreg}).  
Moreover, this definition ensures that the projection onto such a set is well-defined in
its neighbourhood. 
The following remark will be useful later.
\begin{rmrk}
If there exists $\alpha >0$ satisfying $
\xx \in \P_S(\xx+\alpha \vv)$ then $$\forall \beta \geq 0, \  \beta \leq
\alpha, \  \xx \in
\P_S(\xx+\beta \vv).$$ 
\label{pluspetit}
\end{rmrk}
\begin{dfntn}
The proximal subdifferential of function $d_S$ at $\xx $ is the set
$$\partial^P d_S(\xx)= \Big\{\vv \in H,\  \exists M,\  \alpha >0, \
d_S(\yy)-d_S(\xx) +M  |\yy-\xx|^2 \geq
\langle \vv,\yy-\xx \rangle, \ \forall \yy \in B(\xx, \alpha) \Big\} . $$
\end{dfntn}
Let us specify the useful link between the previous subdifferential and the proximal normal cone, which is proved in~\cite{Bounkhel, Clarke2}. 

\begin{prpstn}
The following relation holds true: $$\partial^P d_S(\xx)=N^P(S,\xx) \cap \overline{B(0,1)}.$$
\label{liensousdif}
\end{prpstn}
\begin{rmrk}
 A set $C \subset H$ is convex if and only if it is $\infty$-prox-regular. In this case  $ \NN(C,\xx)=\partial I_{C}(\xx)$  for all $\xx \in C $.
\end{rmrk}

We now come to the main result of this section.
\begin{thrm}
\label{theo:wp}
Assume that $\UU $ is Lipschitz and bounded. Then, for all $T>0$ and all $\qqq_0 \in \ Q $, the
following problem 
$$\left \{ 
\begin{array}{l}
\displaystyle\frac {d \qqq}{dt} + \NNN_\qqq \ni \UU(\qqq) \vsp
\\
\qqq(0)=\qqq_0,
\end{array}
\right.$$
has one and only one absolutely continuous solution $\qqq(\cdot) $ over $[0,T]$.  
\end{thrm}
This well-posedness can be obtained by using results in~\cite{Thibrelax,Thibbv} as soon as we prove that $Q$ is uniformly prox-regular and that the set $\NNN_\qqq$ identifies to the proximal normal cone to $Q$ at $\qqq$. This  is the core of  next subsection.
\begin{rmrk}
 It can be shown that the solution given by Theorem~\ref{theo:wp} satisfies the initial differential equation~(\ref{eqdiff}) (see~\cite{Fred}).
\end{rmrk}

\subsection{Prox-regularity of $Q$}
Let us consider the set  $$ Q_{ij}= \{ \qqq \in  \R^{2N}\virg D_{ij}(\qqq) \geq 0\} .$$
\begin{prpstn}
Let $S$ be a closed subset of $\R^n$ whose boundary $\partial S$ is an oriented $C^2$ hypersurface. For each $\xx \in \partial S$, we denote by $\nu(\xx)$ the outward normal to $S$ at $\xx$.
Then, for each $\xx \in \partial S$, the proximal normal cone to $S$ at $\xx$ is generated by $\nu(\xx)$, i.e. $$ \NN(S,\xx) = \R^+ \nu(\xx).$$ 
\label{conepnlisse} 
\end{prpstn}
\dem
The proof is a straightforward computation  (see~\cite{thesejv}).
\findem

We can also deduce the expression of the proximal normal cone to $Q_{ij}$.
\begin{crllr}
For all $\qqq \in \ Q_{ij} $, 
  $$ \NN(Q_{ij}, \qqq)=-\R^+\GG_{ij}(\qqq) . $$
\end{crllr}
By Definition~\ref{caracpr}, the constant of prox-regularity equals to the largest radius of a ``rolling external ball''. In order to estimate its radius, tools of differential geometry can be used. More precisely, to show that the set $Q_{ij}$ is uniformly prox-regular, we can apply the following theorem, that is proved in~\cite{Delgado}.
\begin{thrm}
Let $C$ be a closed convex subset of $\R^n$ such that $\partial C$ is an oriented $C^2$ hypersurface of $\R^n$. We denote by $\nu_C(\xx)$ the outward normal to $C$ at $\xx$ and by $\rho_1(\xx),..,\rho_{n-1}(\xx) \geq 0 $ the principal curvatures of $C$ at $\xx$.
We suppose that $$ \rho = \sup_{\xx \in \partial C} \ \sup_{ 1 \leq i\leq n-1 } \ \rho_i(\xx) < \infty.$$
Then $S = \R^n \setminus int(C)$ is a $\eta$-prox-regular set with $\eta = \dsp \frac{1}{\rho}$. 
\label{Weingarten}
\end{thrm}

\begin{prpstn}
$Q_{ij}$ is $\eta_0 $-prox-regular with $\eta_0=r\sqrt{2} $.
\label{etaij}
\end{prpstn} 
\dem
The set $int(Q_{ij})$ is obviously the complement of a convex set $C$ which satisfies the assumptions of Theorem~\ref{Weingarten}. The constant of prox-regularity of $Q_{ij}$ can be obtained by calculating its principal curvatures, which are the eigenvalues of Weingarten endomorphism. 
Let $\qqq \in  \partial Q_{ij}$, the outward normal to $C$ at $\qqq$ is equal to $-\nu(\qqq)$, where $$ \nu(\qqq) =-\frac{\GG_{ij}(\qqq)}{\sqrt{2}}=
\frac{\left(0,\dots,0,\ee_{ij}(\qqq),0,\dots,0, -\ee_{ij}(\qqq), 0,\dots,0 \right)}{\sqrt{2}}. $$ Weingarten endomorphism is written as follows,
for all tangent vectors $\hh \in T_{\qqq}(\partial Q_{ij})$, 
$$ \WW_\qqq(\hh):= - \mathrm{D}\nu(\qqq)[\hh]= \dsp \frac{1}{\sqrt{2}|\qq_j-\qq_i|}\left(0,\dots, 0,-\P_{\ee_{ij}^\perp}(\mathrm{h}_j-\mathrm{h}_i),0,\dots,0, \P_{\ee_{ij}^\perp}(\mathrm{h}_j-\mathrm{h}_i), 0, \dots,0 \right) ,$$
with $$\dsp \P_{\ee_{ij}^\perp}(\mathrm{h}_j-\mathrm{h}_i)= (\mathrm{h}_j-\mathrm{h}_i)- [(\mathrm{h}_j-\mathrm{h}_i)\cdot \ee_{ij}
  ]  \ee_{ij}.$$
After some computations, we deduce that the endomorphism $\WW_\qqq $ has two eigenvalues, 0 and $ {\sqrt{2}}/{|\qq_j- \qq_i|} $,
and the latter is equal to 
$1/(r\sqrt 2 )$, which ends the proof.
\findem
Now let us study the set of feasible configurations $Q$, that is the intersection of
 all sets $Q_{ij} $. We begin to determine its proximal normal cone.
 
\begin{figure}
\begin{center}
\includegraphics[width=0.3\textwidth]{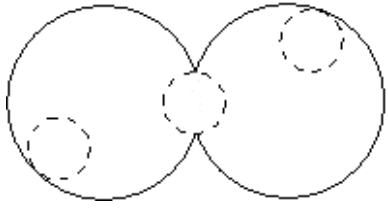}
\caption{Vanishing of the constant of prox-regularity.}
\label{fig:idfausse}
\end{center}
\end{figure}

\begin{figure}
\centering
\psfrag{n1}[l]{$\nn_1$}
\psfrag{n2}[l]{$\nn_2$}
\includegraphics[width=0.2\textwidth]{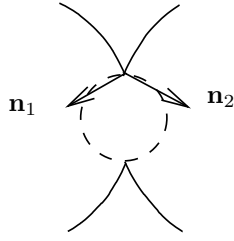}
\caption{Evolution of the angle between vectors $\nn_1$ and $\nn_2$.}
\label{fig:psvectnorm}
\end{figure}

\begin{prpstn}
For all $\qqq \in Q $, $\NN(Q,\qqq )= \sum \NN(Q_{ij},\qqq )= \NNN_{\qqq}
$.
\label{coneprox}
\end{prpstn}
\dem
The second equality follows from~(\ref{exprNq}) and Proposition~\ref{etaij}. Let us prove the first one.
If $ \qqq \in \textmd{int}( Q )$, then for each couple $(i,j)$, $ \qqq \in 
\textmd{int}( Q_{ij} )  $, which implies $$ \NN(Q,\qqq )=
\{0\}=\sum \NN(Q_{ij},\qqq ).$$ We now consider $\qqq \in
\partial Q $ and introduce the following set:
\begin{equation}
 I_{contact}=\{ (i,j), \  i<j,\ D_{ij}(\qqq)=0 
  \}=\{ (i,j), \  i<j,\ \qqq \in \partial Q_{ij} \}.
\label{def:Icontact}
\end{equation}
First, we check that $\NN(Q_{ij},\qqq ) \subset \NN(Q,\qqq ).  $
Let $(i,j) $ belong to $I_{contact} $ (otherwise the previous inclusion is obvious), we consider $\ww \in \NN(Q_{ij},\qqq )\setminus
\{0\}$ and we set  $\vv={\ww}/{|\ww|}$.
By Proposition~\ref{liensousdif},
$\vv \in \partial^P d_{Q_{ij}}(\qqq) $ and thus $$  \exists M,\  \alpha >0, \
d_{Q_{ij}}( \tilde{\qqq})-d_{Q_{ij}}(\qqq) +M |\tilde{\qqq}-\qqq|^2 \geq
\vv \cdot (\tilde{\qqq}-\qqq ), \ \forall \tilde{\qqq} \in B(\qqq,
\alpha). $$ Since $d_{Q_{ij}}(\qqq)=0=d_{Q}(\qqq) $ and
$d_{Q_{ij}}(\tilde{\qqq}) 
 \leq d_{Q}(\tilde{\qqq}) $, it follows that $$  \exists M,\  \alpha >0, \
d_{Q}(\tilde{\qqq})- d_{Q}(\qqq)+M |\tilde{\qqq}-\qqq|^2 \geq
 \vv\cdot (\tilde{\qqq}-\qqq) , \ \forall \tilde{\qqq} \in B(\qqq,
\alpha). $$ Therefore $\vv \in \partial^P  d_{Q}(\qqq) $ and  $\ww \in
\NN(Q,\qqq ).$ Consequently, for each couple $ (i,j) \in
I_{contact} $, we obtain $ \NN(Q_{ij},\qqq )~\subset~\NN(Q,\qqq ) $ as required.
We now want to prove $$ \sum \NN(Q_{ij},\qqq ) \subset \NN(Q,\qqq )  .$$ It suffices to show that $$ \forall \ww_1,\ \ww_2 \in \NN(Q,\qqq )\setminus \{0\} \virg \ww= \ww_1+ \ww_2 \in \NN(Q,\qqq ). $$
Let $\ww_1$ and $\ww_2$ belong to $\NN(Q,\qqq )\setminus \{0\}$, we set $\ww=\ww_1+ \ww_2$, $\vv_1={\ww_1}/{| \ww_1|} $ and $\vv_2={\ww_2}/{|\ww_2|}$.  By Proposition~\ref{liensousdif}, there exists $M_1, M_2 \geq 0$, $\alpha_1, \alpha_2 >0$ such that
$$
\begin{array}{l}
 d_{Q}(\tilde{\qqq})- d_{Q}(\qqq)+M_1 |\tilde{\qqq}-\qqq|^2 \geq
\langle \vv_1,\tilde{\qqq}-\qqq \rangle, \ \forall \tilde{\qqq} \in B(\qqq,
\alpha_1),\vspace{6pt} \\
d_{Q}(\tilde{\qqq})- d_{Q}(\qqq)+M_2 |\tilde{\qqq}-\qqq|^2 \geq
\langle \vv_2,\tilde{\qqq}-\qqq \rangle, \ \forall \tilde{\qqq} \in B(\qqq,
\alpha_2).
\end{array}$$
 So $\ww = |\ww_1| \vv_1+ |\ww_2| \vv_2 $ and 
 the vector
  $ \dsp \vv= {\ww}/{|\ww_1| + |\ww_2|}$ satisfies $|\vv|\leq 1 $. Furthermore $\vv=t\vv_1 +(1-t)\vv_2$, where
 $$\dsp t=\frac{ |\ww_1|}{ (|\ww_1| + |\ww_2|)} .$$  For $\alpha= \min(\alpha_1,\alpha_2) $
 and $M=tM_1 +(1-t)M_2$, the following relation holds
$$d_{Q}(\tilde{\qqq})- d_{Q}(\qqq)+M |\tilde{\qqq}-\qqq|^2 \geq
 \vv \cdot (\tilde{\qqq}-\qqq) , \ \forall \tilde{\qqq} \in B(\qqq,
\alpha). $$ Hence $\vv \in  \partial^P d_{Q}(\qqq) $ and $\ww \in
\NN(Q,\qqq ).$   
To conclude, it remains to check that 
$$ 
\NN(Q,\qqq ) \subset \sum \NN(Q_{ij},\qqq ). 
$$
By~(\ref{decomp}), any $\ww \in \NN(Q,\qqq ) $ can be written $\ww=\vv +\zz=\P_{\NNN_{\qqq}}\ww+\P_{\CCC_{\qqq}}\ww$, with  $\vv \bot \zz$.
Suppose  $\zz
\neq 0$.
Since $\ww \in \NN(Q,\qqq) $, there exists $t>0$ such that $ \qqq \in \P_{Q}(\qqq +t\ww)$. Let
$$s=\min(t,\epsilon) \textmd{   with   }\dsp  \epsilon = \min_{(i,j)
  \notin I_{contact}} \frac{D_{ij}(\qqq)}{\sqrt{2}|\zz|},$$ by Remark~\ref{pluspetit}, we know that $ \qqq \in
  \P_{Q}(\qqq +s \ww)$. Now set
$$\tilde{\qqq}=\qqq + s\ww-s\vv=\qqq+s\zz$$ and show that $ \tilde{\qqq}
  \in Q$. 
 By convexity of 
$D_{ij}$, we have $$ D_{ij}(\tilde{\qqq})\geq D_{ij}(\qqq)+s
  \ \GG_{ij}(\qqq)\cdot \zz, \ \forall (i,j). $$ In addition, for $(i,j)\in
  I_{contact} $, it yields $ \
  \GG_{ij}(\qqq)\cdot \zz \geq 0$, because $\zz \in \CCC_{\qqq} $. Consequently,
  $$\forall (i,j) \in  I_{contact}, \    D_{ij}(\tilde{\qqq})\geq D_{ij}(\qqq)+s
  \ \GG_{ij}(\qqq)\cdot \zz = s
  \ \GG_{ij}(\qqq)\cdot \zz \geq 0    .$$ 
Furthermore, if $(i,j) \notin  I_{contact} $, then $ \dsp s \leq
  \frac{D_{ij}(\qqq)}{\sqrt{2}|\zz|}.$ Hence $$
  D_{ij}(\tilde{\qqq})\geq D_{ij}(\qqq)+s \ \GG_{ij}(\qqq)\cdot \zz \geq
  D_{ij}(\qqq)-s \sqrt{2}|\zz| \geq 0   .$$
That is why $\tilde{\qqq} \in Q$ and $d_{Q}(\qqq +s\ww) \leq |\qqq +s
  \ww - \tilde{\qqq} | = s|\vv |$. Yet $|\qqq +s \ww-\qqq|= s|\ww|> s|\vv|
  $ because $|\ww|^2 =|\vv|^2+|\zz|^2 $. Thus $\qqq \notin
  \P_{Q}(\qqq +s \ww)$, which leads to a contradiction. In conclusion, $\zz=0$ and 
  $\ww =\vv \in \NNN_{\qqq} = \sum \NN(Q_{ij},\qqq ) $, which completes the proof of the proposition.
\findem
Now we want to show the uniform prox-regularity of $Q$. Since $Q$ does not satisfy the same smoothness properties as $Q_{ij}$, the results of differential geometry cannot be applied.
By Theorem~\ref{Weingarten}, if a set is the complement of a smooth convex set, then it is uniformly prox-regular. A natural question arises : Is the intersection of such sets (which is the case for $Q$) uniformly prox-regular with a constant depending only on the constants of prox-regularity of the smooth sets. From a general point of view, this is wrong as illustrated in Figure~\ref{fig:idfausse}. Indeed, we have plotted in solid line the boundary of a set $S$ which is the intersection of two identical disks' complements. This set is uniformly prox-regular but its constant of prox-regularity (equal to the radius of the disk plotted in dashed line) tends to zero when the disks' centres move away from each other. In this situation, the scalar product between normal vectors $\nn_1$ and $\nn_2$ (see Figure~\ref{fig:psvectnorm}) tends to -1. Thus, the constant of prox-regularity of $S$ is also dependent on the angle between vectors $\nn_1$ and $\nn_2$.
We now come to the main result of this subsection: the uniform prox-regularity of $Q$. This result rests on an inverse triangle inequality between vectors $\GG_{ij}(\qqq)$, which is based on angle estimates. 
Let us point out that we do not claim optimality of the constant $\eta$ below.
\begin{prpstn}
$Q$ is $\eta$-prox-regular with 
$$
\eta \sim  \frac {r \sqrt 2 }{2^{3N}} \frac 1 {12^{3N^2}}.
$$
\end{prpstn} 
\dem
 We want to prove (cf.~Proposition~\ref{caracpr}) that there exists $\eta >0$ such that for all 
$\qqq \in Q$ and for all $ \vv \in \NN(Q,\qqq)$,
\begin{equation}
 \vv \cdot (\tilde{\qqq}-\qqq) \leq \frac{|\vv|}{2
  \eta}|\tilde{\qqq}-\qqq|^2, \ \forall \tilde{\qqq} \in Q .
\label{inegfinale}
\end{equation}
By Proposition~\ref{etaij}, for all
$\qqq \in Q_{ij}$ and all $ \ww \in \NN(Q_{ij},\qqq)$, we have
\begin{equation}
 \ww \cdot (\tilde{\qqq}-\qqq)  \leq \frac{|\ww|}{2
  \eta_{0}}|\tilde{\qqq}-\qqq|^2, \ \forall \tilde{\qqq} \in Q_{ij} .
\label{ineginitiale}
\end{equation}
Inegality~(\ref{inegfinale}) is obvious when $\vv =0 $. So we consider $\qqq \in \partial Q$ and $
\vv \in N(Q,\qqq)\setminus \{0\}$.
By Proposition~\ref{coneprox}, $$\vv = -\sum_{(i,j) \in I_{contact}}
\alpha_{ij} \GG_ {ij}(\qqq), \
\alpha_{ij}\geq 0 .$$ 
We recall that $Q \subset Q_{ij}$ so that by~(\ref{ineginitiale}) we obtain
\begin{equation*}
\left( -\sum \alpha_{ij}
\GG_{ij}(\qqq)\right) \cdot (\tilde{\qqq}-\qqq ) \leq \sum\frac{ \alpha_{ij}
  |\GG_{ij}(\qqq)|}{2
  \eta_{0}}|\tilde{\qqq}-\qqq|^2, \ \forall \tilde{\qqq} \in Q.
\end{equation*}
The sum concerns only couples $(i,j)$ belonging to $I_{contact}$ but for convenience, this point is omitted in the notation. As $|\GG_{ij}(\qqq)|= \sqrt{2}$, we get  
\begin{equation*}
\vv \cdot (\tilde{\qqq}-\qqq) \leq 
  \frac{1
  }{\sqrt{2}
  \eta_{0}} \left(\sum \alpha_{ij}\right) |\tilde{\qqq}-\qqq|^2, \
  \forall \tilde{\qqq} \in Q. 
\end{equation*}
To check Inequality~(\ref{inegfinale}), it suffices to find
a constant $\eta >0 $, independent from $\alpha_{ij}$ and from $\qqq$,
satisfying  $$ \left(\sum \alpha_{ij}\right)\frac{1
  }{\sqrt{2} \eta_{0}} \leq  \frac{1}{2 \eta} \left |\sum\alpha_{ij}\GG_{ij}(\qqq)  \right|,
$$ i.e. such that $$\left|\sum \alpha_{ij}\GG_{ij}(\qqq) \right | \geq
\sqrt{2} \frac{\eta}{\eta_{0}} \left(\sum \alpha_{ij}\right) .$$ 
Finally, if we are able to exhibit  $\gamma >0  $ verifying
$$\left|\sum\alpha_{ij}\GG_{ij}(\qqq) \right | \geq
\frac{\sqrt{2}}{\gamma}\left(\sum \alpha_{ij}\right), 
$$ 
then $Q$ will be $\eta$-prox-regular with $$\eta =\dsp\frac{\eta_{0}}{\gamma} = \dsp \frac{r \sqrt{2}}{\gamma}.$$ 
The problem takes the form of  an inverse triangle inequality:
 \begin{equation*} 
\sum\alpha_{ij}| \GG_{ij}(\qqq)  |  = \sqrt{2}\sum \alpha_{ij} 
\leq  \gamma
\left|\sum\alpha_{ij}\GG_{ij}(\qqq)  \right|.
\end{equation*}   
The required result will follow as soon as we prove the main proposition stated below. 
\findem

\begin{prpstn}[Inverse triangle inequality] \mbox{}
\newline There exists $\gamma >1  $ such that for all $\qqq \in Q$,
$$\sum_{(i,j) \in I_{contact}} \alpha_{ij}|\GG_{ij}(\qqq)| 
\leq \gamma \left|\sum_{(i,j) \in I_{contact}} \alpha_{ij}\GG_{ij}(\qqq)
\right|, $$ where $$ I_{contact}=\{(i,j), \ i<j,\  D_{ij}(\qqq)=0 \} 
\textmd{ and } \alpha_{ij} \textmd{ are nonnegative reals}. $$
Constant $\gamma $ can be fixed as follows  
$$\gamma = 
\dsp \left[
\frac{1}{2}\left(
1- \left (
1+\left(\dsp \frac{1}{12^{2N}}\right)
\right )   ^{-1/2}
\right) 
\right]^{\dsp -\frac{3N}{2}} .$$


\label{inegtrianginv}
\end{prpstn}

\begin{rmrk}
 Note the sign of coefficients $\alpha_{ij}$. From a general point of view, this inequality is obviously wrong if  these coefficients are just assumed real. Indeed, for $N$ large enough, the cardinal of the set $I_{contact}$ could be strictly larger than  $2N$, which induces a relation between vectors $\GG_{ij}(\qqq)$ (see Fig.~\ref{fig:crist} for such a degenerate situation). 
\end{rmrk}
The following elementary lemma asserts an inverse triangle inequality for two vectors.
\begin{lmm}
Let $u_1 $ and $u_2 $ be two vectors of $\R^{2N}$ satisfying $u_1
\cdot u_2= \cos \theta |u_1||u_2|,$ with $\cos \theta >-1 $. Then for all $$ \nu \geq  \nu_{\theta}:=\sqrt{\frac{2}{1+\cos \theta}} $$ we have $|u_1| +|u_2| \leq \nu |u_1+u_2|. $
\label{inegtrianginv3}
\end{lmm}
Proof of the inverse triangle inequality:
We propose here a method based on angle estimates with vectors $\GG_{ij}(\qqq)$ as pointed out in Figure ~\ref{fig:psvectnorm}. We use a recursive proof on the number of involved vectors. We are going to check that there exists $\delta >1$ such that for all subset $I \subset I_{contact} $ and for all $\alpha_{ij} >0$,
 $$\sum_{(i,j) \in I \subset I_{contact}} \alpha_{ij}|\GG_{ij}(\qqq)| 
\leq \delta^{|I|} \left|\sum_{(i,j) \in I \subset I_{contact}}
\alpha_{ij}\GG_{ij}(\qqq) 
\right|. $$ 
\emph{Initialization:} Suppose that the cardinality of $I$ equals to 1, in other words,
$I =\{(i,j)\}$. So we clearly have for all $\alpha_{ij}>0$ and all $\delta >1$, 
\begin{equation}
\alpha_{ij} |\GG_{ij}(\qqq)|   
= |\alpha_{ij} \GG_{ij}(\qqq)| \leq \delta
|\alpha_{ij}\GG_{ij}(\qqq)| .  
\label{eq:rec} 
\end{equation}
\emph{Recursion assumption:} \\ If $|J|=p$, then we have for all $\alpha_{ij} >0$
\begin{equation} 
\sum_{(i,j) \in J \subset I_{contact}} \alpha_{ij}|\GG_{ij}(\qqq)| 
\leq \delta^p \left|\sum_{(i,j) \in J \subset I_{contact}} \alpha_{ij}\GG_{ij}(\qqq)
\right|. 
\label{eq:hyprec}
\end{equation}
Take a subset $ I \subset I_{contact}$ with $|I|=p+1$. For any $$\ww= \sum_{(i,j) \in I} \alpha_{ij} \GG_{ij}(\qqq) ,
$$  
 with $\alpha_{ij} >0$, we choose $(k,l)  \in
I$ and define $J= I \setminus~\{(k,l) \} $, $$\ww_1= \sum_{(i,j)\in J} \alpha_{ij} \GG_{ij}(\qqq) \textmd{ and } \ww_2=
\alpha_{kl} \GG_{kl}(\qqq). $$
We need the following lemma which will be later proved.
\begin{lmm} If $\ww_1 \neq 0$, the following inequality holds
$$\dsp \frac{\ww_1 \cdot \ww_2}{|\ww_1 ||\ww_2 |
} \geq - \kappa, \textmd{ with } \kappa= \left (1+ \left( \dsp \frac{1}{12}\right)^{2N} \right ) ^{-1/2}. $$ 
\label{prodscal}
\end{lmm}
Consequently, if $\ww_1 \neq 0$, from Lemma~\ref{inegtrianginv3}, we deduce $| \ww_1| + |\ww_2 | \leq \dsp  \sqrt{\frac{2}{1-\kappa}}|
 \ww_1  + \ww_2 |$ (this inequality obviously holds for $\ww_1=0$). By denoting  $\delta= \dsp \sqrt{\frac{2}{1-\kappa}} >1$, we get
\begin{equation}
| \ww_1 |  + | \ww_2 | \leq \delta |\ww|.
\label{ineglem}
\end{equation}
Applying recursion assumption (\ref{eq:hyprec}) and (\ref{ineglem}), we obtain
 $$ \sum_{(i,j) \in
  I_{contact}} 
\alpha_{ij}|\GG_{ij}(\qqq)|\leq \alpha_{kl} |\GG_{kl}(\qqq) |
+\delta^p  | \ww_1|\leq \delta^p \left( | \ww_2 |  +  |
\ww_1|\right) \leq \delta^{p+1} | \ww|, $$  which ends the proof of (\ref{eq:rec}) by recursion. 
 As $|I_{contact}| \leq \dsp 3N $, the inverse triangle inequality is checked with 
$\gamma =\delta^{3N}$. 
\findem
Proof of Lemma~\ref{prodscal}: 
It suffices to deal with $\ww_2 = \GG_{kl}(\qqq) .$
By setting $$\beta_{ij}=\left\{\begin{array}{ll}
\alpha_{ij}& \textmd{ if } i<j \vsp \\
\alpha_{ji}& \textmd{ else, }
\end{array}\right.$$ we have  $$ \ww_1=(F_1,F_2,...,F_N) \textmd{
  where } F_p= \sum \beta_{ip} \ee_{ip}.$$ Thus, $F_k \in \R^2 $ 
can be interpreted as a pressure force exerted on the $k^{\textmd{th}}$ person by its neighbours (different from the individual $l$). Similarly, $-F_k$ can be seen as a reaction force.
We are looking for a lower bound of
  $$\Delta_{kl}:= \dsp  \frac{\ww_1\cdot \ww_2}{|\ww_1 | |\ww_2 |} =  \dsp \frac{- F_k
  \cdot \ee_{kl}- F_l \cdot \ee_{lk} }{\sqrt{2} \sqrt{  \sum_{i=1}^N
  |F_i|^2}}  .$$
\begin{figure}
\begin{subfigure}[Case 1]{
\psfrag{e}[l]{$\ee_{lk}$}
\psfrag{f}[l]{$\ee_{kl}$}
\psfrag{g}[l]{$-F_k$}
\psfrag{h}[l]{$-F_l$}
\psfrag{k}[l]{$q_l$}
\psfrag{j}[l]{$q_k$}
\includegraphics[scale=1]{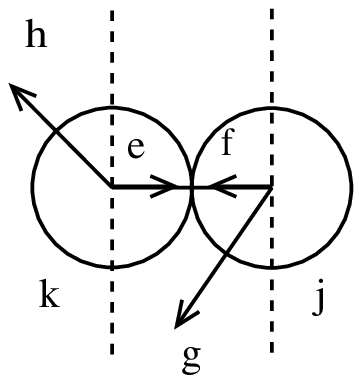}
\label{fig:cass1}}
\end{subfigure}
\begin{subfigure}[Case 2a]{
\psfrag{e}[l]{$\ee_{kl}$}
\psfrag{f}[l]{$\ee_{lk}$}
\psfrag{g}[l]{$-F_l$}
\psfrag{h}[l]{$-F_k$}
\psfrag{k}[l]{$q_k$}
\psfrag{j}[l]{$q_l$}
\includegraphics[scale=1]{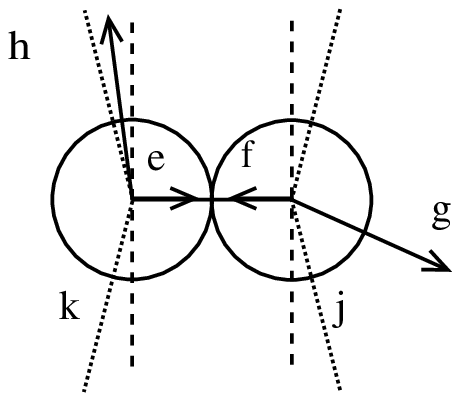}
\label{fig:cas2a}}
\end{subfigure}
\begin{subfigure}[Case 2b]{
\psfrag{e}[l]{$\ee_{kl}$}
\psfrag{f}[l]{$\ee_{lk}$}
\psfrag{g}[l]{$-F_l$}
\psfrag{h}[l]{$-F_k$}
\psfrag{k}[l]{$q_k$}
\psfrag{j}[l]{$q_l$}
\includegraphics[scale=1]{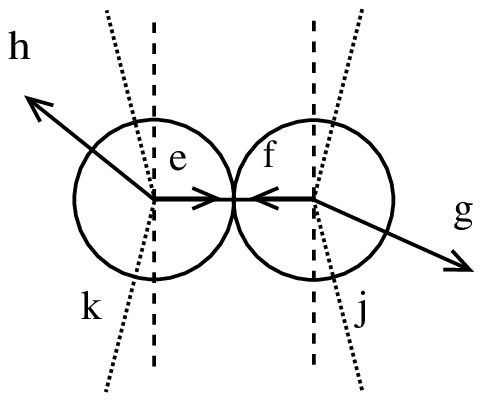}
\label{fig:cas2b}}
\end{subfigure}
\end{figure}
\emph{Case 1:} $- F_k \cdot \ee_{kl} \geq 0 $ or $-
  F_l \cdot \ee_{lk} \geq 0$  \\ \\ Suppose that, for example (cf figure \ref{fig:cass1}) $ - F_k \cdot \ee_{kl} \geq 0 $.
   Using $|F_l \cdot
  \ee_{lk} | \leq |F_l|,$  we get  
  $$\Delta_{kl} \geq \frac{-F_l \cdot
  \ee_{lk}}{\sqrt{2}\sqrt{\sum |F_i|^2}}   \geq \frac{-1}{\sqrt{2}}. $$
In this case, $ \dsp \kappa = 2^{-1/2}.$ \\
\emph{Case 2:} $- F_k \cdot \ee_{kl} <0 $ and $-
  F_l \cdot \ee_{lk} < 0$  \\
 \emph{Case 2a:}  $- F_k \cdot \ee_{kl} \geq \dsp
  -\frac{1}{4} |F_k| $ or $- F_l \cdot \ee_{lk} \geq \dsp
  -\frac{1}{4} |F_l| $ \\
 Suppose that, for example (cf Figure~\ref{fig:cas2a}), $ \dsp - F_k \cdot \ee_{kl} \geq \dsp
  -\frac{1}{4} |F_k| $.
It can be shown that $$-\frac{1}{4} \leq \frac{- F_k \cdot
  \ee_{kl}}{\sqrt{\sum |F_i|^2}} \textmd{   and   } \frac{- F_l \cdot \ee_{lk}}{\sqrt{\sum |F_i|^2}} \geq -1,$$ which yields 
  $$\Delta_{kl} \geq  
  \frac{1}{\sqrt{2}}\left(-\frac{1}{4}-1 \right)= -\frac{5}{4
  \sqrt{2}}>-1. $$
In this case $ \dsp \kappa={5}/({4\sqrt{2}})  $.
 \\ \\
 \emph{Case 2b:}  $- F_k \cdot \ee_{kl} < \dsp
  -\frac{1}{4} |F_k| $ and $- F_l \cdot \ee_{lk} < \dsp
  -\frac{1}{4} |F_l| $ (cf Figure~\ref{fig:cas2b}).\\
We need the following lemma.
\begin{lmm}
There exists $\tilde{k} $ and $\tilde{l} $
different from $k$ and $l$ verifying $\tilde{k} \neq \tilde{l} $ and 
$$\begin{array}{lll}
\dsp |F_{\tilde{k}}|& \geq &\epsilon |F_k|,\vsp \\
\dsp  |F_{\tilde{l}}|& \geq &\epsilon |F_l|,
\end{array}$$ with $\epsilon =  {1}/{12}^{2N}. $ 
\label{lembarbare}
\end{lmm}
 We deduce that $$\sum |F_i |^2 \geq \dsp |F_k|^2+
|F_l|^2+|F_{\tilde{k}}|^2+ 
|F_{\tilde{l}}|^2\geq (1+\epsilon^2) \left[|F_k|^2+
|F_l|^2 \right]. 
$$
 Therefore $$ |\Delta_{kl}| \leq \dsp
\frac{1}{\sqrt{1+\epsilon^2}} 
\left(\frac{|F_k|+|F_l|  }{\sqrt{2} \sqrt{|F_k|^2+
|F_l|^2}  } \right)\leq \frac{1}{\sqrt{1+\epsilon^2}}  . $$
In this case, $ \dsp \kappa=\frac{1}{\sqrt{1+\epsilon^2}},$ which concludes the proof of Lemma~\ref{prodscal}.
\findem

\begin{figure}
\psfrag{e}[l]{$\ee_{kl}$}
\psfrag{f}[l]{$\ee_{lk}$}
\psfrag{g}[l]{$-F_l$}
\psfrag{k}[l]{$q_k$}
\psfrag{j}[l]{$q_l$}
\psfrag{a}[l]{$q_{k_1}$}
\psfrag{b}[l]{$q_{k_2}$}
\psfrag{c}[l]{$q_{k_i}$}
\begin{center}
\includegraphics[scale=1]{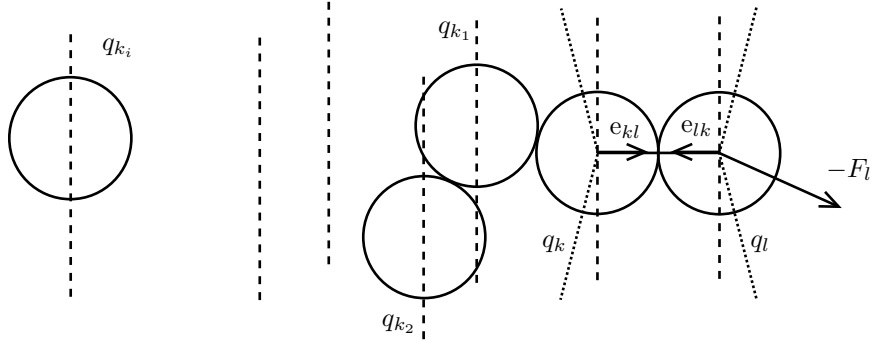}
\caption{Construction of sequence $(k_i) $}
\label{fig:suite}
\end{center}
\end{figure}

Proof of Lemma~\ref{lembarbare}:
We firstly consider  $$ -F_k = \sum_{i=1}^{V_k}
\beta_{kj_{0,i}} \ee_{kj_{0,i}},$$ where $V_k$ is the number of neighbours of individual $k$ (individual $l$ excepted) ($V_k\leq 5$). As a consequence, $$- F_k \cdot
\ee_{kl} = \sum_{i=1}^{V_k} 
\beta_{kj_{0,i}} \ee_{kj_{0,i}} \cdot \ee_{kl} . $$  There exists $ k_1
\in\{j_{0,1}, j_{0,2},...,j_{0,V_k}  \}$ ($ k_1 \neq k,l$) such that for all $i \in \{ 1,..., V_k\}$
$ \beta_{kk_1} \ee_{kk_1}\cdot \ee_{kl} \leq \beta_{kj_{0,i}}
\ee_{kj_{0,i}}\cdot \ee_{kl} . $ It is obvious that 
$$ \beta_{kk_1} \ee_{kk_1}\cdot \ee_{kl}< -\frac{1}{6}F_k \cdot
\ee_{kl}\leq -\frac{1}{24} |F_k|.$$  In fact, individual $ k_1$ is the neighbour who exerts the largest pressure force on person $k $.  As illustrated in Figure~\ref{fig:suite}, individual $k$ is between persons $l$ and $k_1$.\\
If $\dsp |F_{k_1}| \geq \frac{1}{48} |F_k| $, then we set
$\tilde{k}=k_1. $ Else
 $\dsp |F_{k_1}|< \frac{1}{48} |F_k| $, and we produce the same reasoning with $$-F_{k_1}= \beta_{k_1k}
\ee_{k_1k} + \sum_{i=1}^{V_{k_1}}
\beta_{k_1j_{1,i}} \ee_{k_1j_{1,i}}, $$ where $V_{k_1} \leq 5$. Thus, $$-F_{k_1} \cdot \ee_{kl}= \beta_{k_1k}
\ee_{k_1k}\cdot \ee_{kl} + \sum_{i=1}^{V_{k_1}}
\beta_{k_1j_{1,i}} \ee_{k_1j_{1,i}}\cdot \ee_{kl} .$$ Since $ \dsp -\beta_{k_1k}
\ee_{k_1k}\cdot \ee_{kl} < -\frac{1}{24}|F_k| $ and $ \dsp -F_{k_1}\cdot \ee_{kl} 
\leq |F_{k_1}|< \frac{1}{48} |F_k| $, we obtain $$\sum_{i=1}^{V_{k_1}}
\beta_{k_1j_{1,i}} \ee_{k_1j_{1,i}}\cdot \ee_{kl}=-F_{k_1} \cdot \ee_{kl}
-\beta_{k_1k} 
\ee_{k_1k}\cdot \ee_{kl} < -\frac{1}{48}|F_k| .$$
As previously, there exists $ k_2
\in\{j_{1,1}, j_{1,2},...,j_{1,V_{k_1}}  \} $ ($k_2 \notin \{k,k_1 \}
$), such that $$ \beta_{k_1k_2} 
\ee_{k_1k_2}\cdot \ee_{kl} < -\frac{1}{4 \times 12^2}|F_k|$$ 
(Similarly, see Figure~\ref{fig:suite}, individual $k_1$ is between persons $k_2$ and $k$). \\
If $|F_{k_2}| \geq  \dsp \frac{1}{4} \left(\frac{1}{12} \right)^2 |F_k|$,
 we set $\tilde{k}=k_2 .$
 Else, we continue by defining a sequence $(k_i)$ (cf Figure~\ref{fig:suite}) such that
$$\left \{\begin{array}{l}
\dsp k_0=k \vsp \\
\dsp |F_{k_{i+1}}| < \frac{1}{4} \left(\frac{1}{12} \right)^{i+1} |F_k|\vsp \\
\dsp \beta_{k_ik_{i+1}} 
\ee_{k_ik_{i+1}}\cdot \ee_{kl} < -\frac{1}{4} \left(\frac{1}{12}
\right)^{i}\frac{1}{6} |F_k|.
\end{array}
\right. $$
It can be shown that $\dsp k_{i+1} \notin \{k_0,
k_1,..k_i \} . $ This construction ends at most in $N-2$ steps:
$$\exists m<N-1 \textmd{ satisfying } |F_{k_m}| \geq \frac{1}{4}
\left(\frac{1}{12} \right)^m |F_k|. $$ Finally we set $$\tilde{k}=k_m .$$
Analoguously, we deal with $F_l$, by constructing a sequence $(l_i)$ verifying similar properties. We can check that $ \tilde{k}\neq \tilde{l} $ in proving that $$\{k_0,k_1,..k_m\} \cap\{l_0,l_1,..l_p\} = \emptyset .$$  
The proof of Lemma~\ref{lembarbare} is achieved by taking $
\epsilon ={1}/{12}^N$.
\findem
\section{Numerical scheme}
\label{sec:numscheme}
\subsection{Time-discretization scheme}
We present in this section a numerical scheme to approximate the solution
 to~(\ref{incldiff}). 
 The numerical scheme we propose is based on a first order expansion of
the constraints expressed in terms of velocities.  
The time interval is denoted by $[0,\ T] $. Let $ N\in \mathbb{N}^\star$, $h = T/N $ be the time
 step and $t^\itt=\itt h $ be
 the computational times. 
 We denote by 
$\qqq^\itt $ the approximation of  $\qqq(t^\itt )$. 
The  next  configuration is obtained as
$$
\qqq^{\itt +1} = \qqq^\itt  + h \ \uu^{\itt },
$$
where $$\uu^{\itt }=  \P_{\CCC^h_{\qqq^\itt }} (\UU(\qqq^\itt )) \hbox{ with} $$
$$\CCC^h_{\qqq} =  \{ \vv \in \R^{2N} \virg D_{ij}(\qqq) + h \ \GG_{ij}(\qqq)\cdot \vv  \geq  0 \quad \forall \, i<j \}. $$
The scheme can be also interpreted in the following way. Let us introduce the set 
$$
\Qc(\qqq)  = \{ \tilde \qqq \in \R^{2N} \virg D_{ij}(\qqq) + \GG_{ij}(\qqq)\cdot
(\tilde \qqq - \qqq) \geq 0 \quad \forall \, i<j 
\},
$$
which can be seen as an inner convex approximation of $Q$ with respect
to $\qqq$. 
Note that $\Qc(\qqq)$ is defined in such a way that $Q$ is the union
of all sets $\Qc(\qqq) $, $\qqq \in Q $. The scheme can be expressed in terms of position:
$$ \qqq^{\itt +1} = \P_{\Qc(\qqq^\itt )} (\qqq^\itt  + h \UU(\qqq^\itt )) .$$ 
In this form it appears as a prediction-correction algorithm: predicted position vector $\qqq^\itt  + h \UU(\qqq^\itt )$, that may not be admissible,  is projected  onto the approximate set of feasible configurations.  

\begin{rmrk}
It is straightforward to check that
\begin{equation}
\label{eq:diffincdisc}
\frac {\qqq^{\itt +1} - \qqq^\itt }{h} + \NN(\Qc(\qqq^\itt ),\qqq^{\itt +1} ) \ni
\UU(\qqq^\itt ),
\end{equation}
so that the scheme can also be seen as a 
semi-implicit discretization  of~(\ref{incldiff}),
 where $ \NN(\Qc(\qqq^\itt ),\qqq^{\itt +1} ) $ approximates $ \NN(Q,\qqq^{\itt } )$.
 \end{rmrk}

 Convergence  of this scheme shall be proven in a 
 forthcoming paper.
\subsection{Numerical solutions}
\begin{figure}
\centering
\includegraphics[width=0.25\textwidth]{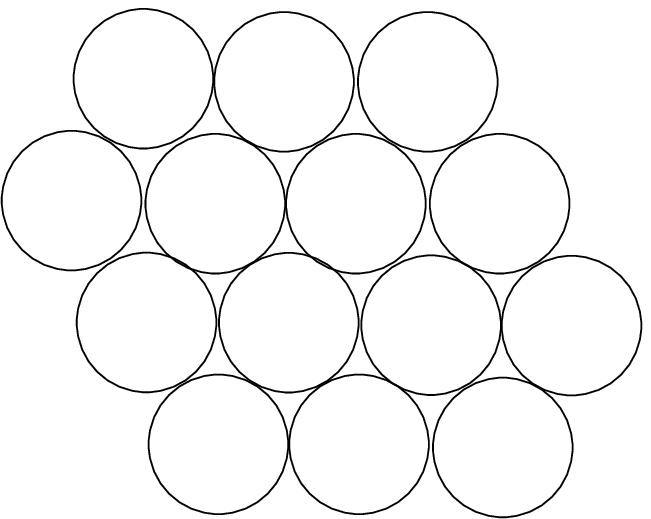} 
\caption{A case of non-uniqueness for Kuhn-Tucker multipliers.}
\label{fig:crist}
\end{figure}

In the model, the discrete actual velocity $\uu^\itt $ is the projection of the spontaneous
velocity onto the approximated set of feasible velocities. We propose here to
solve this projection by a Uzawa algorithm (note that any algorithm could
be used to perform this task). For convenience,  explicit dependence of vectors and matrices upon the current configuration is omitted (e.g. $\UU$ stands for $\UU(\qqq^n)$, $D_{ij}$ for $D_{ij}(\qqq^n)$, etc\dots).
The actual velocity $\uu $ solves the following minimization problem under constraints 
  $$ \dsp \uu=
 \underset{\vv \in \CCC^h_{\qqq}}{\mathrm{argmin}}\  |\vv - \UU|^2 .$$ 
 Uzawa algorithm is based on a
reformulation of this minimization problem in a saddle-point form.
We introduce the associated Lagrangian
$$
  L \left(\vv,\bfmu \right) = \frac{1}{2}
|\vv-\UU|^2 -
    \sum_{1 \leq i<j \leq N} \mu_{ij}\  \left(D_{ij} + h\  \GG_{ij}\cdot  \vv \right ). $$
and the following linear mapping
$$
\begin{array}{llll}
    \BB  : & \mathbb{R}^{2N} & \rightarrow &  \mathbb{R}^{\frac{N(N-1)}{2}}
    \\
    & \vv & \mapsto & -h \left(\GG_{ij}\cdot   \vv \right)_{i<j}
\end{array}
 $$
%
With these notations, the set $\CCC^h_{\qqq}$ can be written:
\begin{eqnarray*}
   \CCC^h_{\qqq}&=&\left\{\vv \in \mathbb{R}^{2N}\ ,\ \forall \bfmu \in \left(
    \mathbb{R}^+\right)^{\frac{N(N-1)}{2}}\ , \
    - \sum_{1 \leq i<j \leq N} \mu_{ij} \left( \ D_{ij} + h\  \GG_{ij} \cdot \vv \right)\leq 0\right\}
    \\
     &=& \left\{\vv \in\mathbb{R}^{2N}\ ,\ \forall \bfmu \in \left(
    \mathbb{R}^+\right)^{\frac{N(N-1)}{2}}\ , \
     \bfmu \cdot (\BB \vv-\Dvec ) 
    \leq 0\right\}.
\end{eqnarray*}
where $\Dvec = \Dvec (\qqq)\in  \mathbb{R}^{N(N-1)/2}$ is the vector of distances.
The existence of a saddle-point 
$$(\uu, \bflambda) \in \mathbb{R}^{2N} \times (\R^+)^{\frac{N(N-1)}{2}}
$$
 for this problem is well-known (see e.g.~\cite{Ciarlet}) and it is characterized by the next system:
$$\left\{ \begin{array}{l}
           \uu  +\phantom{}^{t}\BB \bflambda = \UU  \vsp \\
           \bfmu \cdot( \BB\uu- \Dvec )   \leq 0 \virg \forall \bfmu  \geq 0\vsp \\
           \bflambda  \cdot ( \BB\uu - \Dvec )  = 0.
          \end{array}
 \right. $$
Uzawa algorithm produces two sequences $(\vv^k)\in \left(\R^{2N}\right)^\mathbb{N}$ and $(\bfmu^k)\in \left((\R^+)^{\frac{N(N-1)}{2}}\right)^\mathbb{N} $ according to
\begin{eqnarray*}
    \bfmu^0 & = & 0
    \\
    \vv^{k+1} & = & \UU - \phantom{}^{t}\BB\bfmu^k
    \\
    \bfmu^{k+1} & = & \Pi_+\left(\bfmu^k + \rho \left[
    \BB\vv^{k+1}- \Dvec  \right] \right),
\end{eqnarray*}
where $\Pi_+$ is the euclidean projection onto the cone of vectors with nonnegative
components (a simple cut-off in practice), and $ \rho> 0 $ is a fixed parameter.
The algorithm can be shown to converge as soon as $0 < \rho < 2/\|\BB\|^2$ (see~\cite{Ciarlet}).
More precisely, the sequence $(\vv^k)$ converges to $\uu$ and it can be shown that the sequence $(\bfmu^k)$ tends to some $\overline{\bflambda} \in~(\R^+)^{\frac{N(N-1)}{2}} $ such that $(\uu, \overline{\bflambda})$ is a saddle-point of $L$. Notice that in general, the Kuhn-Tucker multiplier $\bflambda $ is not unique as illustrated in Figure~\ref{fig:crist}. In this case, the configuration of 14 people shows 29 contacts, consequently matrix  $\phantom{}^{t}\BB$ is not injective. 

\begin{rmrk}{(Link between local prox-regularity and speed of convergence for Uzawa algorithm)}  We denote by $\Gmat $ the matrix whose columns are vectors $\GG_{ij}$, where $(i,j) \in I_{contact}$ (defined by ~(\ref{def:Icontact})), and we introduce 
$\CC= \phantom{}^{t}{\Gmat } \Gmat$.
The size of this square matrix is equal to $n_{contact}$ which is the cardinal of $I_{contact}$.
By inverse triangle inequality (see Proposition~\ref{inegtrianginv}),  there exists a constant $\gamma$ such that 
for all $\bflambda \in (\R^+)^{n_{contact}}$ satisfying $|\bflambda|_1=1 $, we have $$ \left| \sum \lambda_{ij} \GG_{ij}\right|^2=\phantom{}^{t}{\bflambda} \phantom{}^{t}{\Gmat } \Gmat  \bflambda =\phantom{}^{t}{\bflambda} \CC\bflambda \geq \frac{2}{\gamma^2}. 
$$
We define,  for $\qqq \in Q$, a local parameter $\gamma_\qqq$ satisfying $$\min_{\genfrac{}{}{0pt}{}{|\bflambda|_1 = 1}{ \bflambda \geq 0} }
 \phantom{}^{t}{\bflambda} \CC\bflambda = \dsp \frac{2}{\gamma_\qqq^2} ,
 $$
  and  $\eta_\qqq = r\sqrt{2}/\gamma_\qqq $.
 Let us show that parameter $\eta_\qqq$ (setting a lower bound of the local prox-regularity of $Q$ at point $\qqq $) and the condition number of matrix $\CC$  are closely related when $\CC$ is non-singular.
By denoting $\eta_{min}$ the smallest eigenvalue of $\CC$, it follows that $$\eta_{min} = \min_{|\bflambda|_2=1} \phantom{}^{t}{\bflambda} \CC\bflambda =  \min_{|\bflambda|_2 \geq 1} \phantom{}^{t}{\bflambda} \CC\bflambda \leq \min_{\genfrac{}{}{0pt}{}{|\bflambda|_2 \geq 1}{ \bflambda \geq 0} }
 \phantom{}^{t}{\bflambda} \CC\bflambda.$$
Since for all $ \bflambda$, $|\bflambda|_1 \leq  \sqrt{n_{contact}} |\bflambda|_2  $, we have $$\min_{\genfrac{}{}{0pt}{}{|\bflambda|_2 \geq 1}{ \bflambda \geq 0} }
 \phantom{}^{t}{\bflambda} \CC\bflambda \leq \min_{\genfrac{}{}{0pt}{}{|\bflambda |_1 \geq \sqrt{\mathit{n_{contact}}}}{ \bflambda \geq 0} }
 \phantom{}^{t}{\bflambda} \CC\bflambda = \mathit{n_{contact}} \min_{\genfrac{}{}{0pt}{}{|\bflambda|_1 \geq 1}{ \bflambda \geq 0} }
 \phantom{}^{t}{\bflambda} \CC\bflambda .$$
Finally, $$\eta_{min} \leq n_{contact} \min_{\genfrac{}{}{0pt}{}{|\bflambda|_1 \geq 1}{ \bflambda \geq 0} }
 \phantom{}^{t}{\bflambda} \CC \bflambda = \mathit{n_{contact}} \min_{\genfrac{}{}{0pt}{}{|\bflambda|_1 = 1}{ \bflambda \geq 0} }
 \phantom{}^{t}{\bflambda} \CC\bflambda = \dsp \frac{2 \mathit{n_{contact}} }{\gamma_\qqq^2} .$$
Thus $$ \eta_{min} \leq \dsp \frac{ 6N }{\gamma_\qqq^2}.$$ 
Furthermore, the condition number of matrix $\CC $ equals to $$ cond_2 (\CC) =\| \CC \|_2 \| \CC^{-1} \|_2 =\frac{\eta_{max}}{\eta_{min}}. 
$$ 
 Since $|\GG_{ij}(\qqq)|= \sqrt{2}$, we obtain $\| \CC \|_2 = \eta_{max} \geq 2 $, 
hence $$ cond_2 (\CC)\geq  \dsp \frac{2}{\eta_{min}} \geq \frac{2 \gamma_\qqq^2}{\mathit{6 N}}  
\geq
 \frac{4r^2}{ 6\eta_\qqq^2 { N} },
 $$
 which quantifies how the condition number of $A$ varies with  $\eta_\qqq$.
Since the matrix appearing in Uzawa algorithm is  $\CC=\phantom{}^{t}{\Gmat } \Gmat  $, we expect that this algorithm converges less quickly for configurations with low local prox-regularity. In numerical simulations, we noticed indeed that solving the saddle-point problem requires more iterations in case of a jam.
\label{conditionnement}
\end{rmrk}

\section{Numerical results}
\label{sec:resnum}

In order to illustrate the contact model, we propose here an example of spontaneous velocity. The choice of the spontaneous velocity is important because this velocity reflects pedestrian behaviour. A lot of choices are obviously possible. The spontaneous velocity of an individual has to take into account obstacles in the room and specify how he wants to get around them.
So this velocity depends on the room's geometry but it can be made dependent on other people positions too. Indeed, it is possible here to integrate individual strategies (deceleration or jam's avoiding). We refer the reader to~\cite{esaim, MauryvenelTGF, justrat} for other examples of spontaneous velocity. Here we restrict ourselves to simple behavourial model: people tend to optimize their own path, regardless of others. 

\subsection*{An example of spontaneous velocity}
We consider here the simplest choice for the spontaneous velocity. All the individuals have the
same behaviour: they want to reach the exit by following the shortest
path avoiding obstacles. Then, the spontaneous velocity's expression can be specified:
$$ 
\dsp
 \UU(\qq)=( \UU_0(\qq_1),\ldots, 
 \UU_0(\qq_N))
\textmd{  with  } \UU_0(\xx)= -  \textmd{s} \ \nabla \mathcal{D}(\xx), 
 $$
where $\mathcal{D}(\xx)$
represents the geodesic distance between the position $\xx$ and the nearest
exit and $\textmd{s} >0$ denotes the
speed. 
\begin{figure}
\centering
 \includegraphics[width=0.46\textwidth]{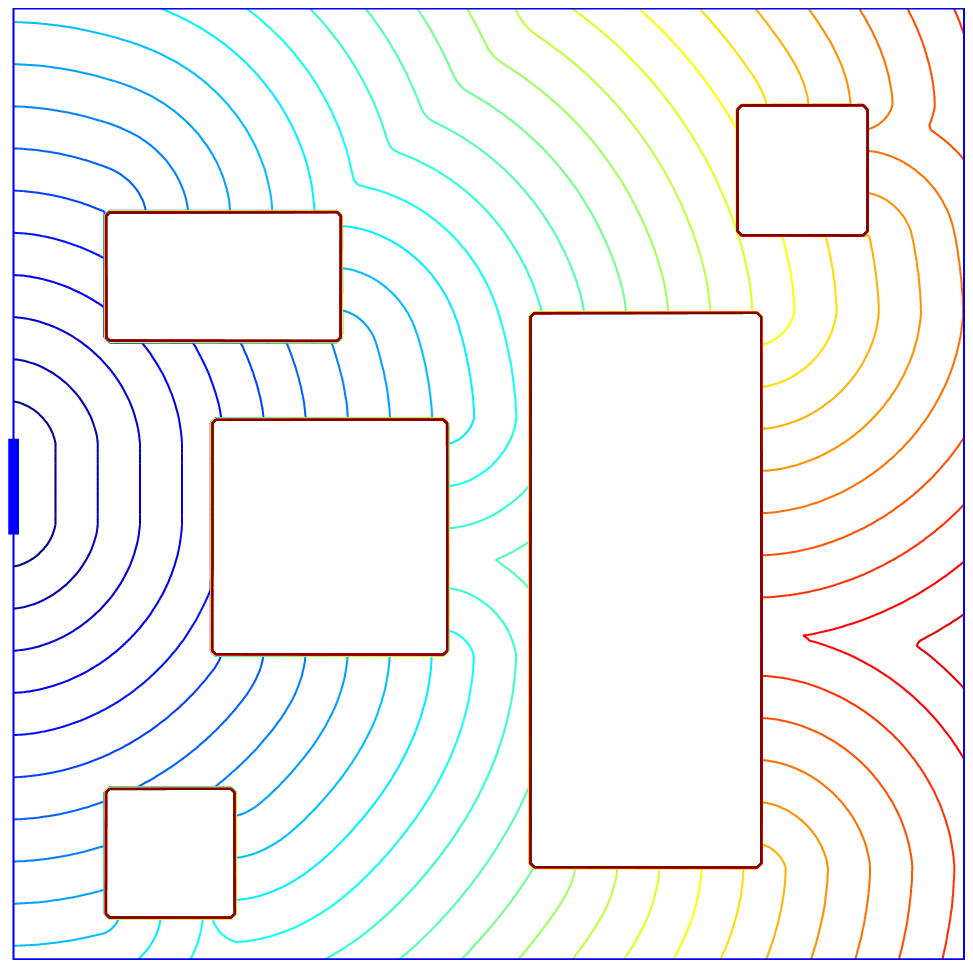}
 \includegraphics[width=0.475\textwidth]{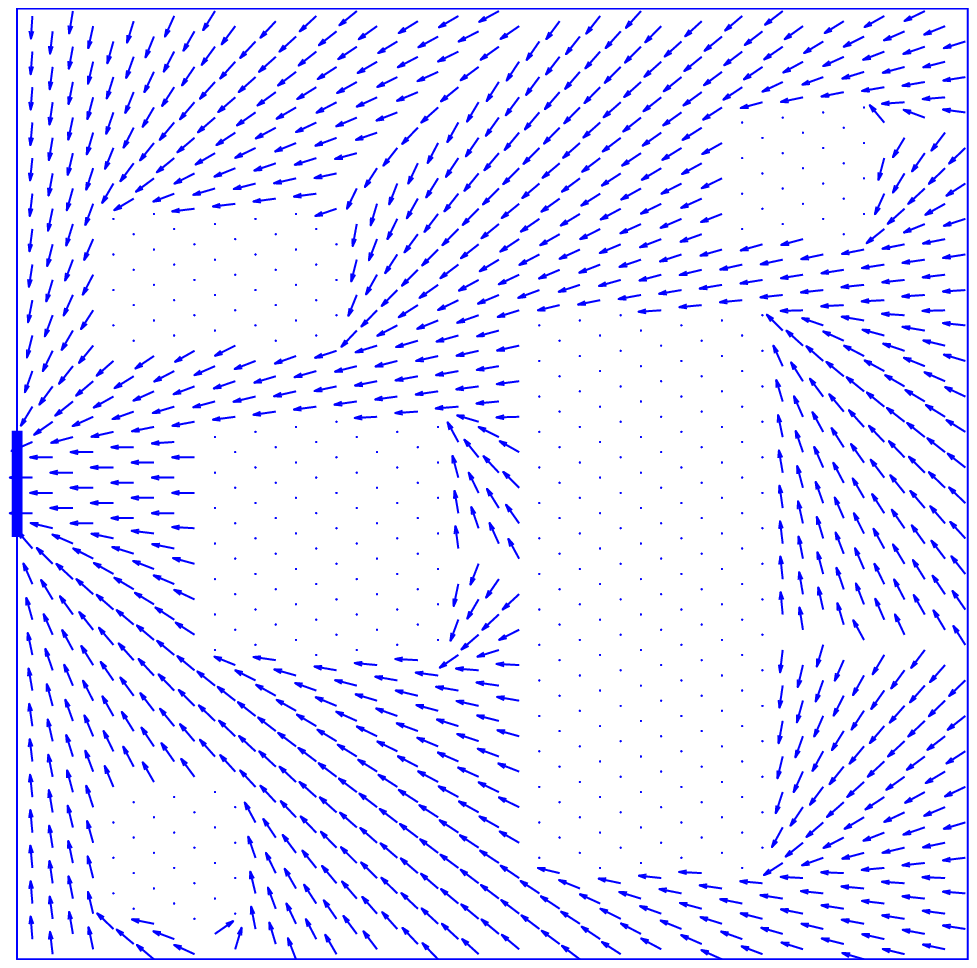}
\caption{Contour levels of the geodesic distance $\mathcal{D}$ and velocity field $\UU_0$.}
\label{fig:contlev}
\end{figure} 
In order to compute $\mathcal{D}$, we have used the Fast
Marching Method introduced by R. Kimmel and J. Sethian 
in \cite{FFM}.
In this method, the value of $\mathcal{D}$ is computed at each point
of a grid. The value at the exit's nodes is set to zero. Then,
the values of the distance at the other points is computed step by
step so that a discrete version of  $|\nabla \mathcal{D}
 | =1$ is satisfied. Moreover, the distance at the nodes situated in
 the obstacles is fixed to a large value, which prevents the shortest
 path from going across them.
In Figure~\ref{fig:contlev}, we have considered a room with 5 obstacles and the exit is situated to the left.
We note that by following the built velocity field, people are going to avoid obstacles. 

Our aim is to simulate evacuation of any building consisting of several floors. We have chosen an object oriented programming method and we have implemented this Fast Marching Method in a C++ code. Let us detail this code. On each floor, the spontaneous velocity is directed by the shortest path avoiding obstacles to the nearest exit or stairwell. In the stairs, people just want to go down.
We have integrated this spontaneous velocity in the C++ code SCoPI: Simulations of Collections of Interacting Particles developped by A. Lefebvre (see~\cite{Aline, Aline_part}). This code allows us to compute the actual velocity as the projection of the spontaneous velocity as described in Section~\ref{sec:numscheme}.

\begin{rmrk}
Notice that the velocity field produced by this strategy is not continuous as soon as the room is not convex, which  rules out  Theorem~\ref{theo:wp}. This lack of regularity  is not important in practical applications~: the places at  which  it occurs (in particular upstream obstacles) are  emptied after a few moments. The main consequence is the discontinuity of the future configurations  with respect to initial data, which is not surprising from a modelling standpoint.
\end{rmrk}

We propose to illustrate the behaviour of the algorithm in two situations. 
The first one corresponds to a many-individual evacuation from a square  room through a single exit, the second one illustrate the capability of the approach to handle complicated geometries. 
For these two experiments, it will be noticed that the contacts between the individuals and the obstacles have to be handled (as the contacts between people). Even if an individual want to avoid an obstacle, he can be pushed on it by people behind them. 

\subsection*{Simple evacuation}
\begin{figure}
\begin{center}
\begin{tabular}{c@{\hspace{1cm}}c}
\resizebox{!}{!}{\includegraphics[width=0.32\textwidth]{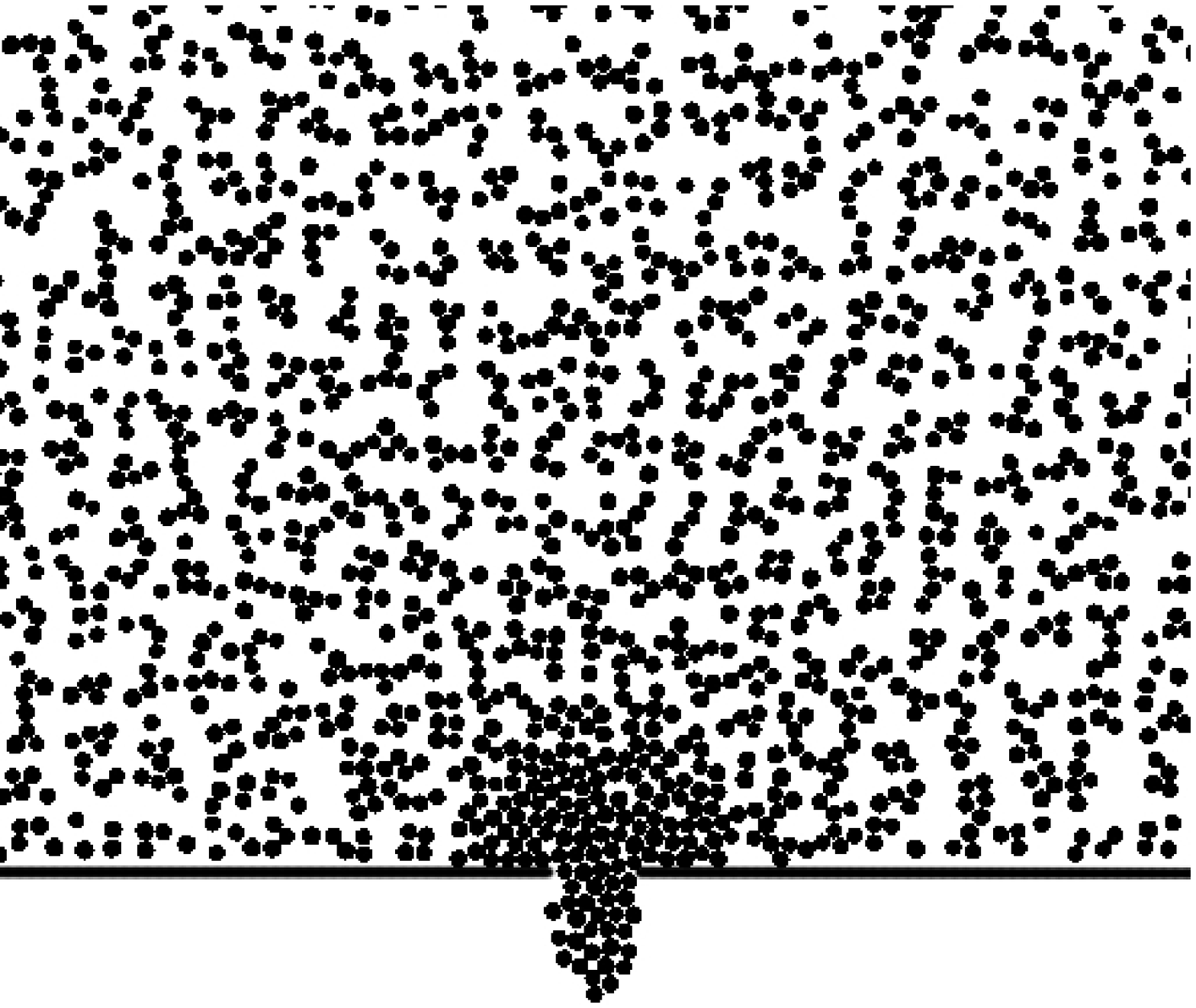}}
&
\resizebox{!}{!}{\includegraphics[width=0.32\textwidth]{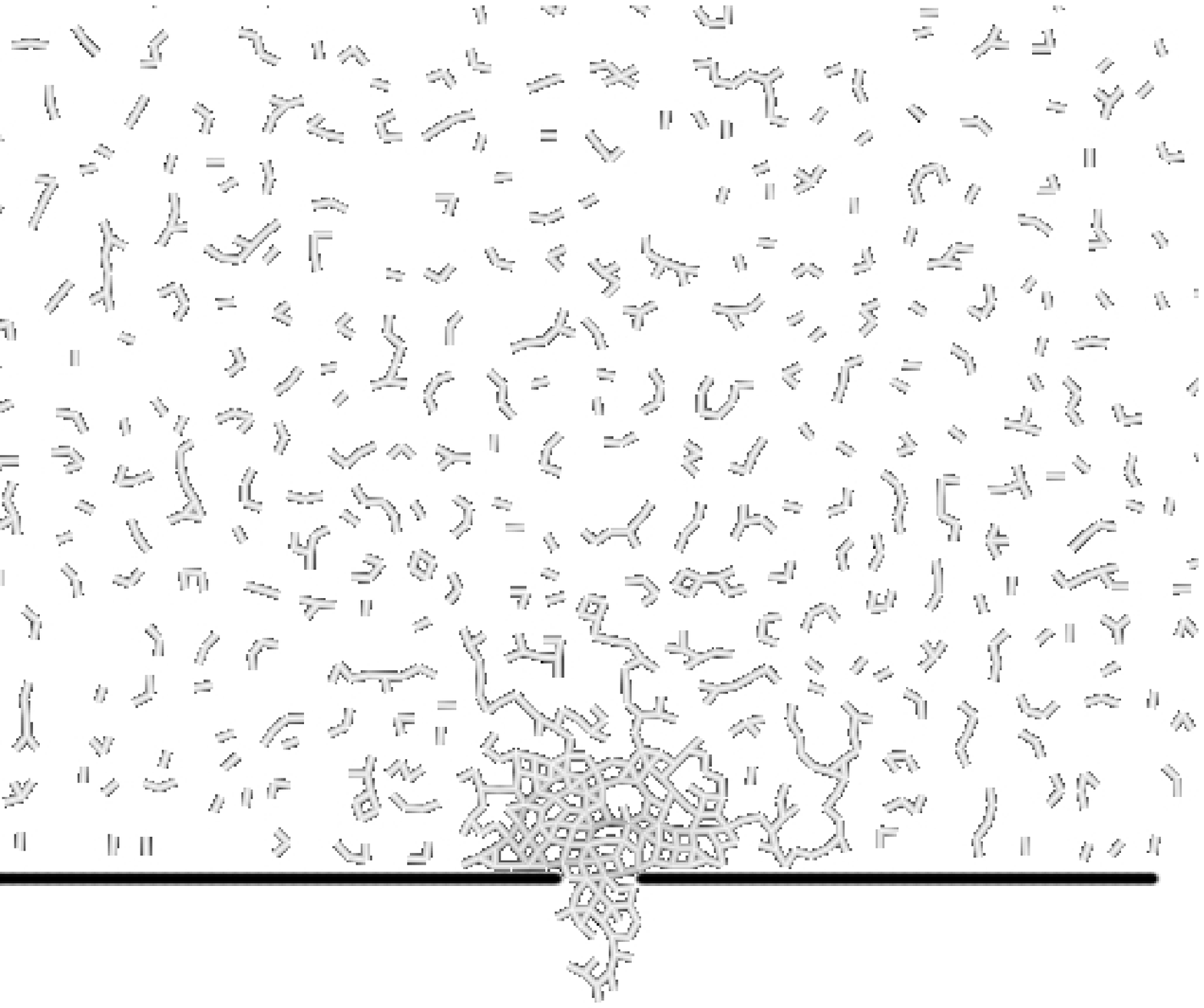}}
\end{tabular}
\\
\begin{tabular}{c@{\hspace{1cm}}c}
\resizebox{!}{!}{\includegraphics[width=0.32\textwidth]{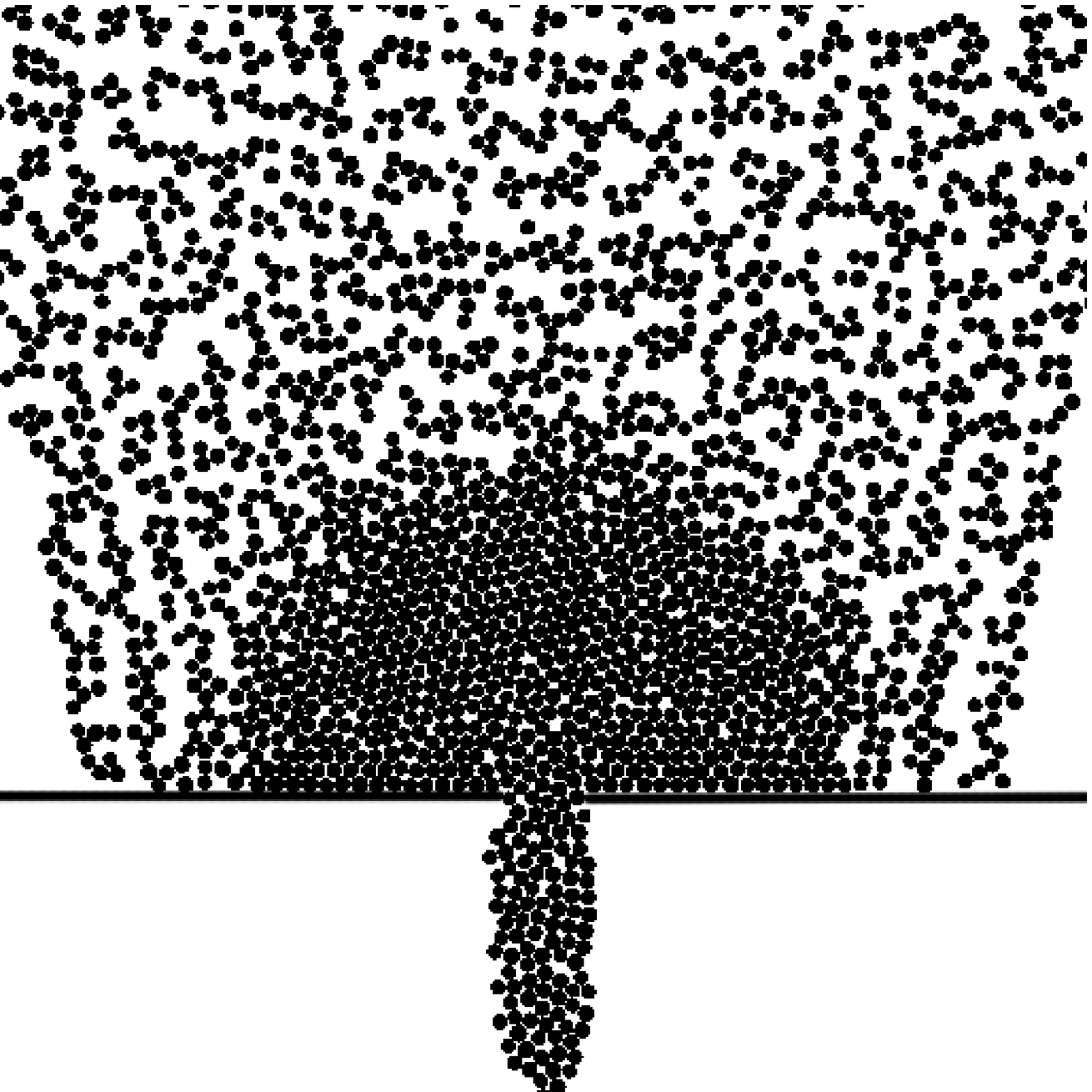}}
&
\resizebox{!}{!}{\includegraphics[width=0.32\textwidth]{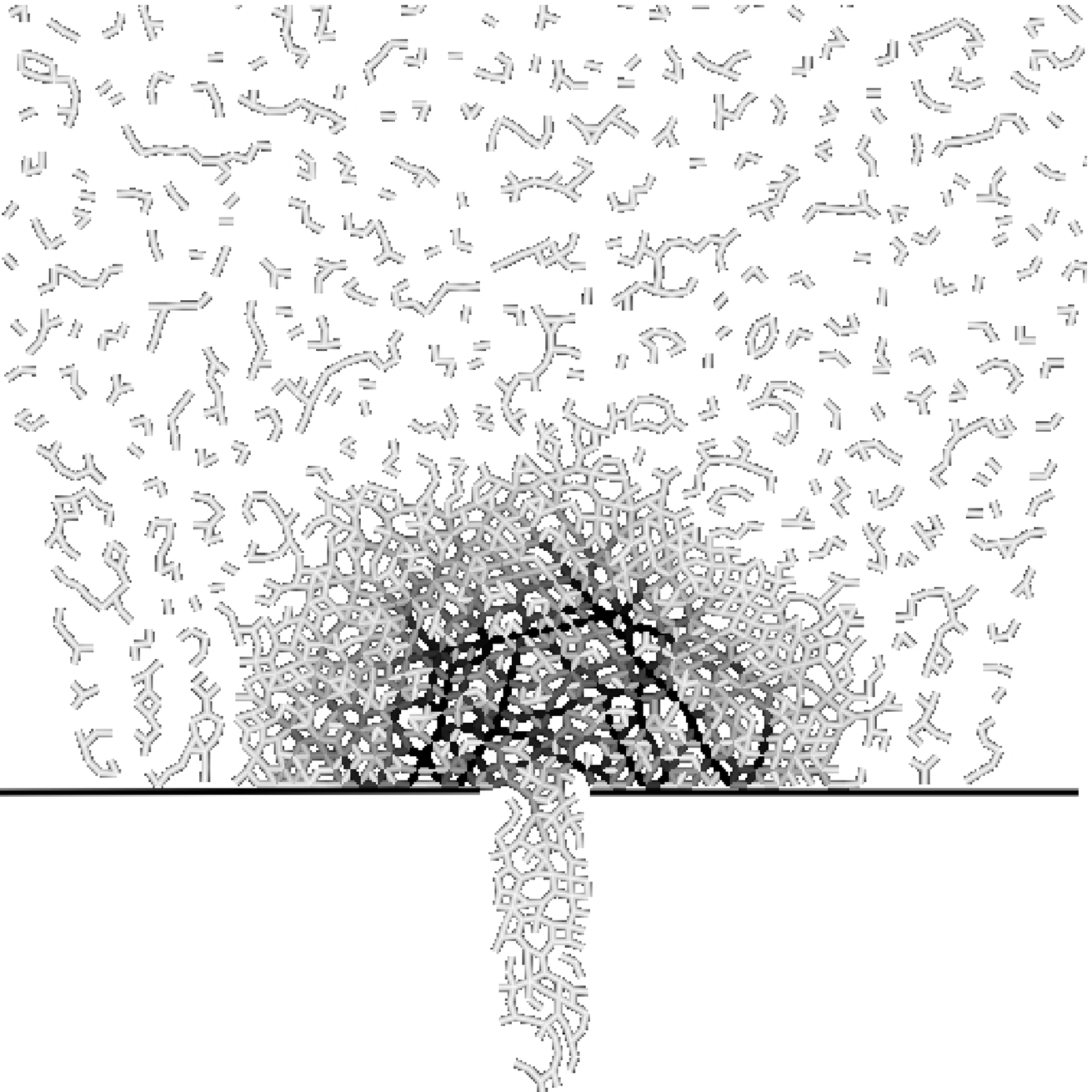}}
\end{tabular}
\\
\begin{tabular}{c@{\hspace{1cm}}c}
\resizebox{!}{!}{\includegraphics[width=0.32\textwidth]{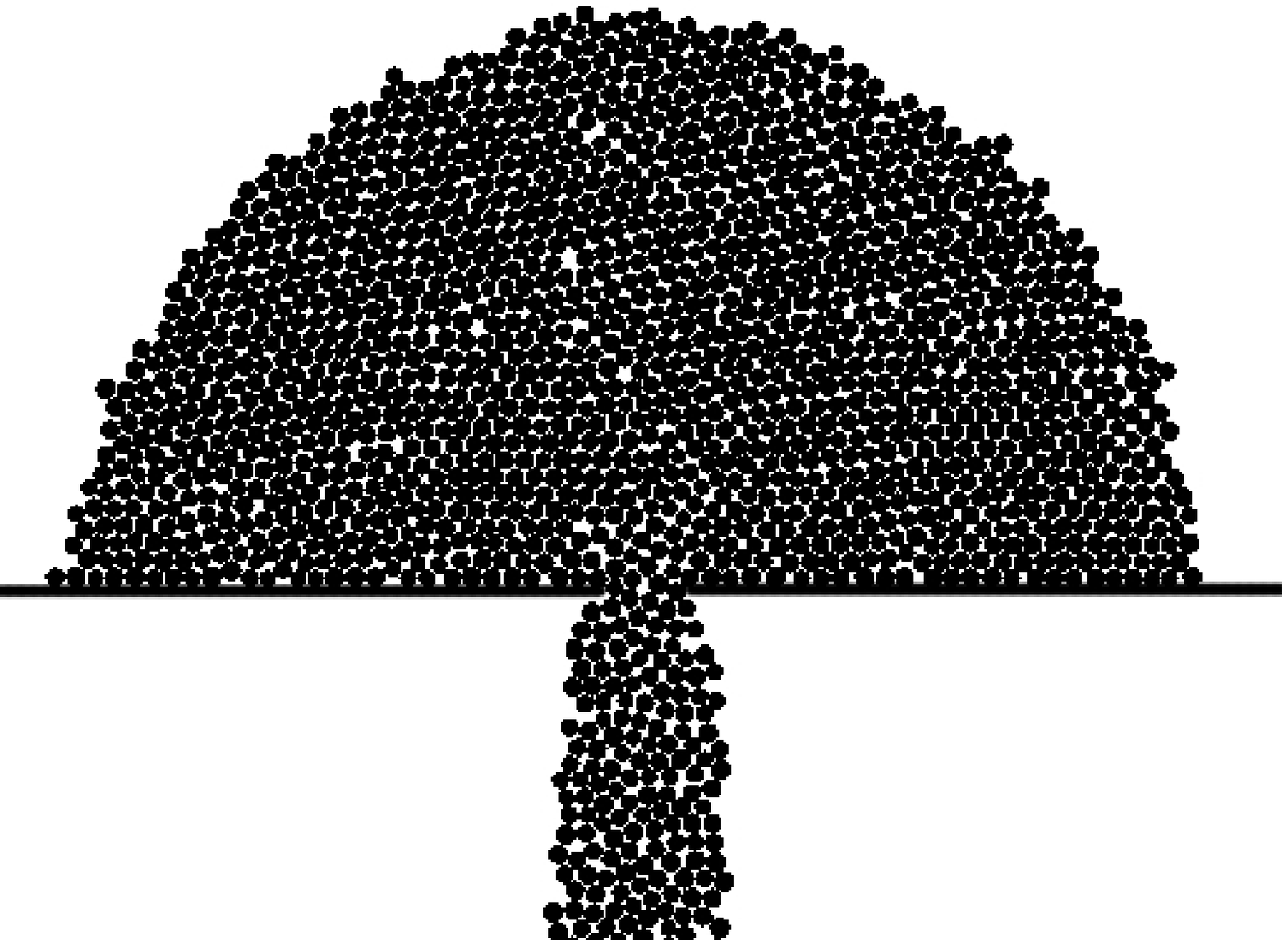}}
&
\resizebox{!}{!}{\includegraphics[width=0.32\textwidth]{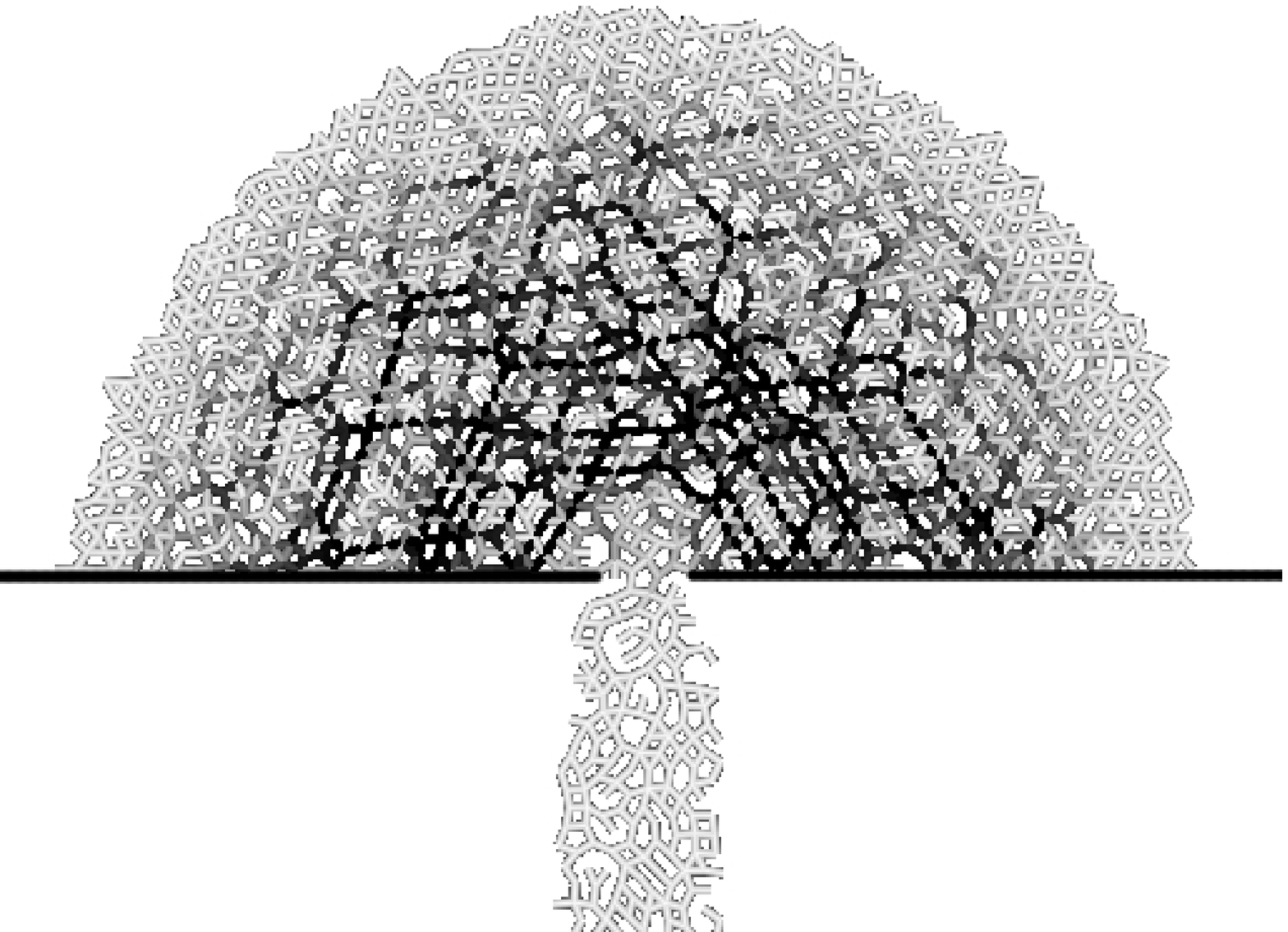}}
\end{tabular}
\\
\begin{tabular}{c@{\hspace{1cm}}c}
\resizebox{!}{!}{\includegraphics[width=0.32\textwidth]{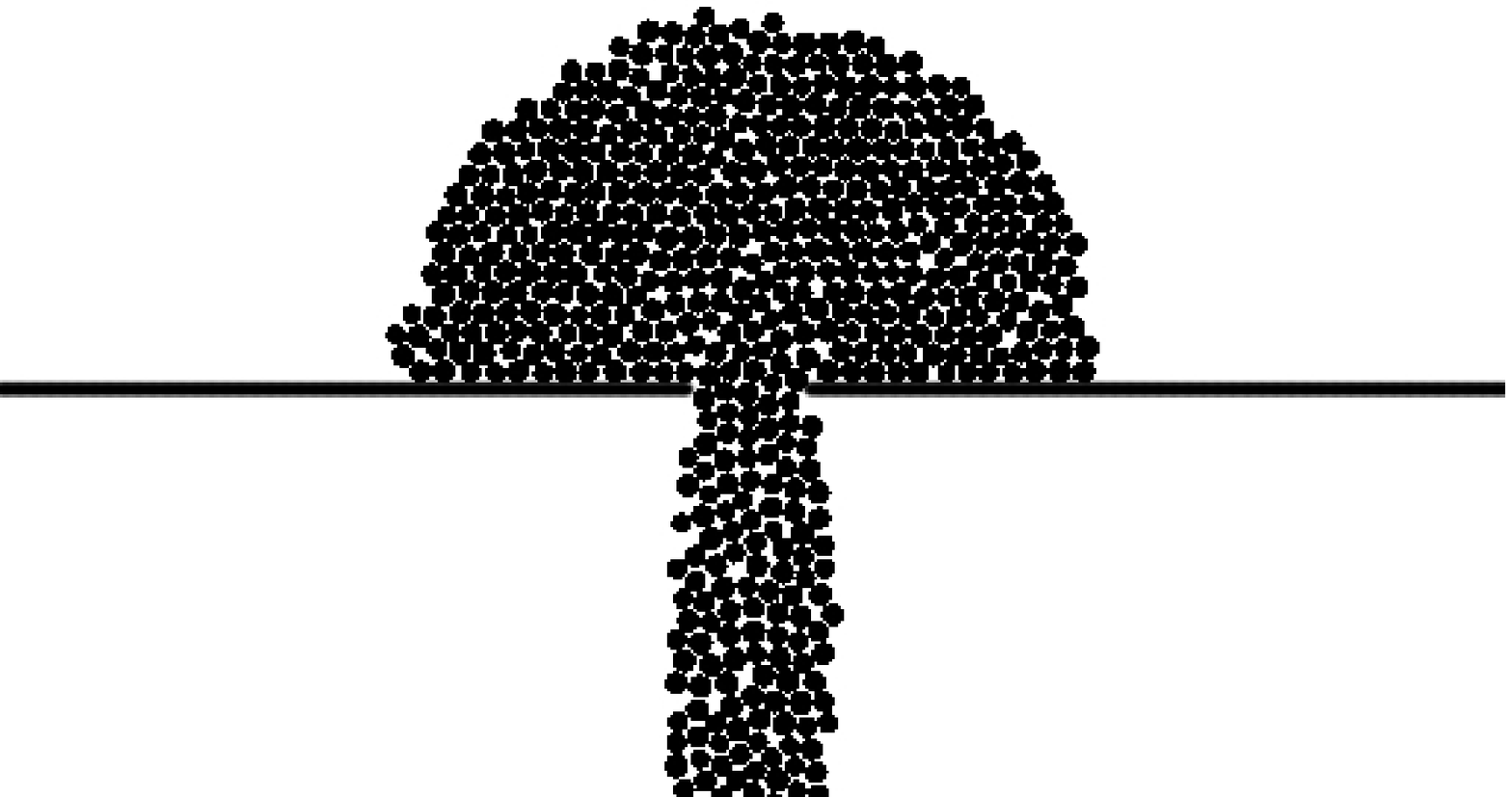}}
&
\resizebox{!}{!}{\includegraphics[width=0.32\textwidth]{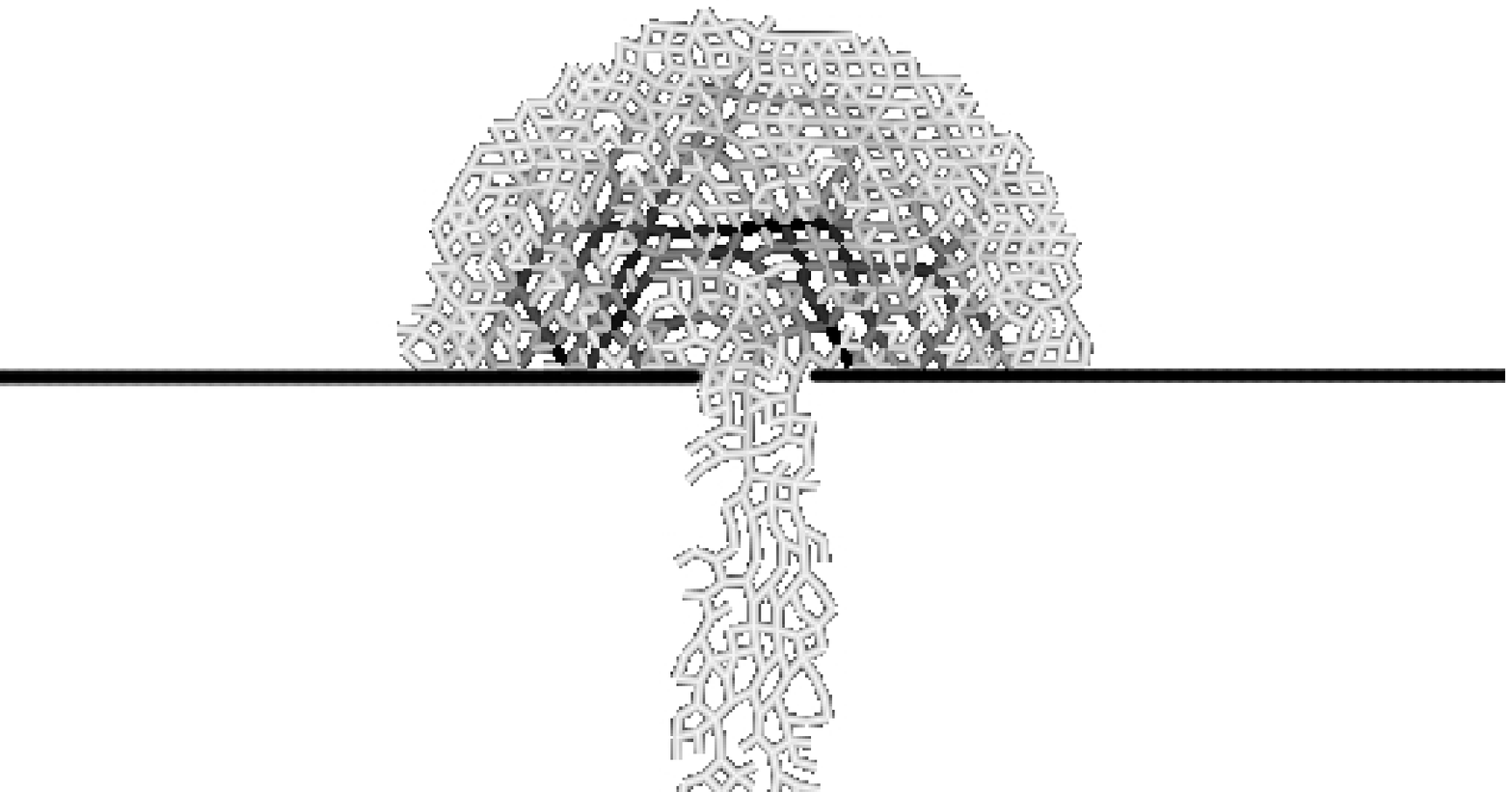}}
\end{tabular}
\end{center}
\caption{Arches.}
\label{fig:arcades}
\end{figure}

We consider the situation of 1000 people which are randomly distributed over a square room. The spontaneous velocity field corresponds straight pathlines towards the exit at constant speed. As the field has a negative divergence, it tends to increase the local density, so that congestion is rapidly reached in the neighbourhood of the exit, and the congestion front propagated upstream as long as it is feeded by incoming people. 
In Figure~\ref{fig:arcades}, we represented the current configuration and the corresponding network of interaction pressures: for any couple of disks
in contact, we represent the segment between centers, having its color (from
white to black) depend upon the (positive) Kuhn-Tucker multiplier which handles
the corresponding contraint. We recover the apparition of arches upstream the exit.
The Kuhn-Tucker multipliers $\lambda_{ij}$ quantify the way $\UU$, the spontaneous velocity
field, does not fit the constraints, and as such they can be interpreted in terms
of pressures undergone by individuals. Although it would be presumptuous at
this stage to assimilate $\lambda_{ij}$ to an actual measure of the discomfort experienced
by persons $i$ and $j$, it is obvious that high values for those Kuhn-Tucker multipliers
can be expected on zones where people are likely to be crushed.

\subsection*{Complex geometry}
\begin{figure}
\centering
\includegraphics[width=0.4\textwidth]{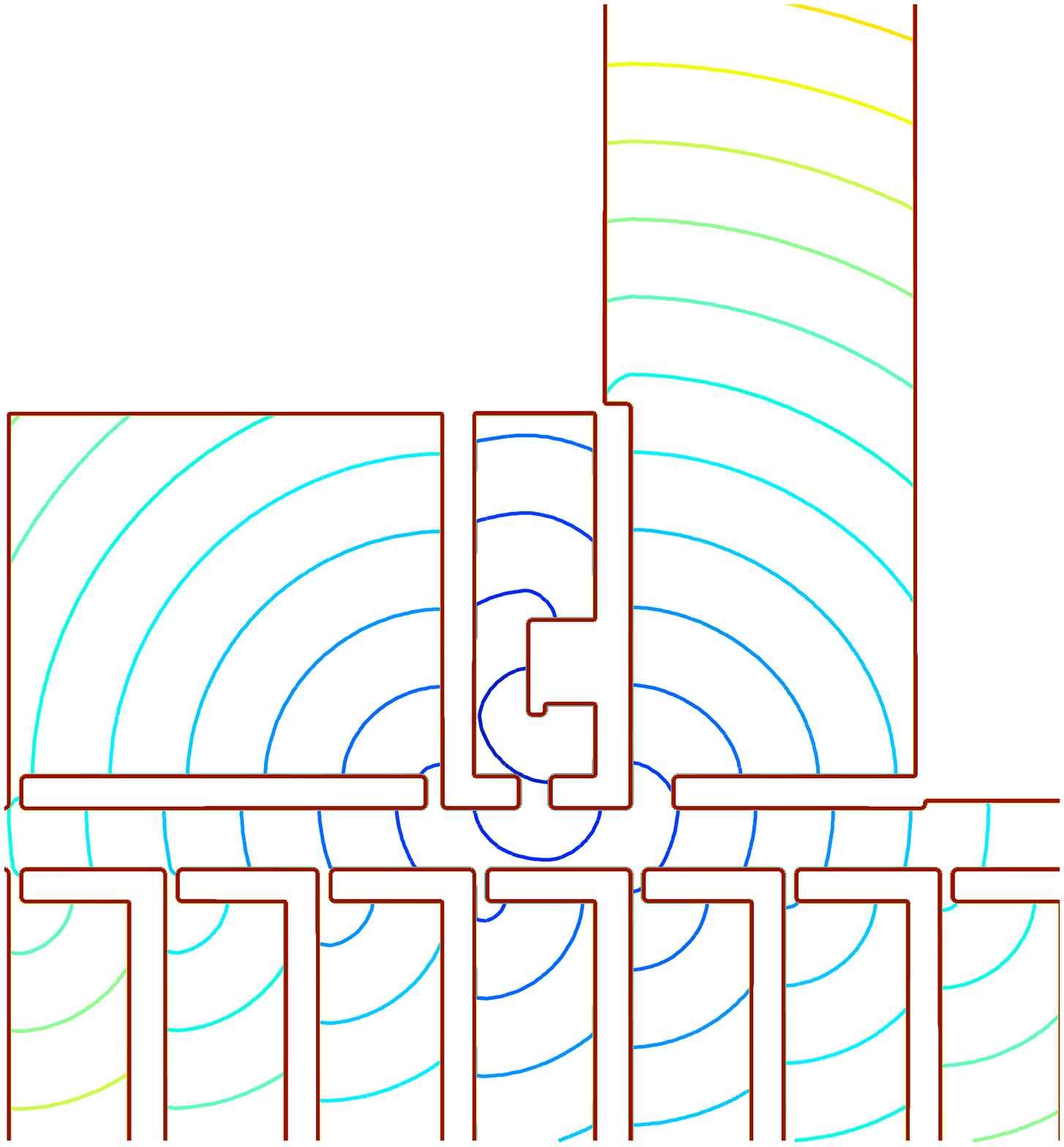}
\caption{Geometry and isovalues for the  geodesic distance.}
\label{fig:batmath}
\end{figure}

\begin{figure}
\begin{center}
\begin{tabular}{c@{\hspace{1cm}}c}
\resizebox{!}{!}{\includegraphics[width=0.4\textwidth]{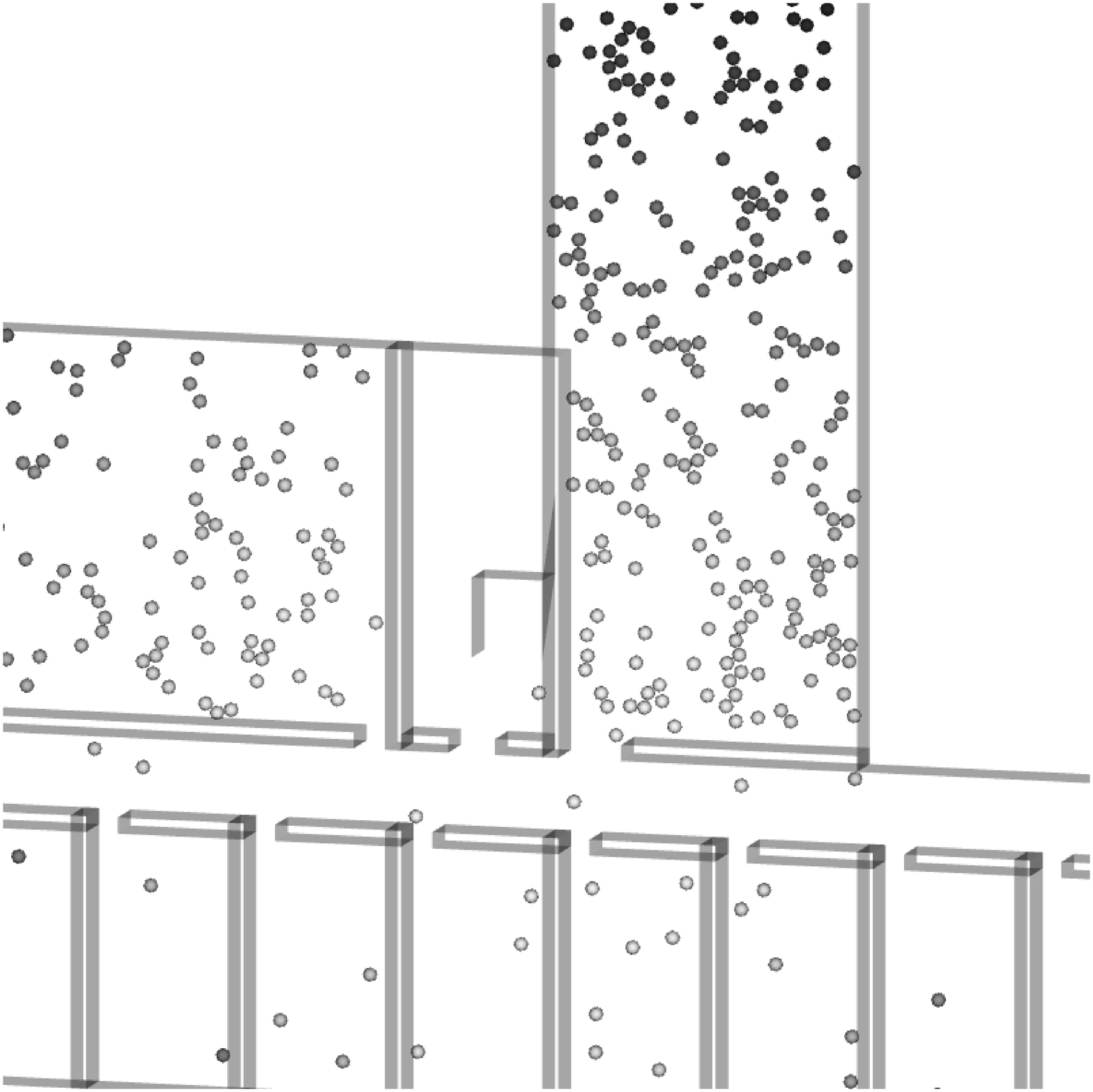}}
& 
\resizebox{!}{!}{\includegraphics[width=0.4\textwidth]{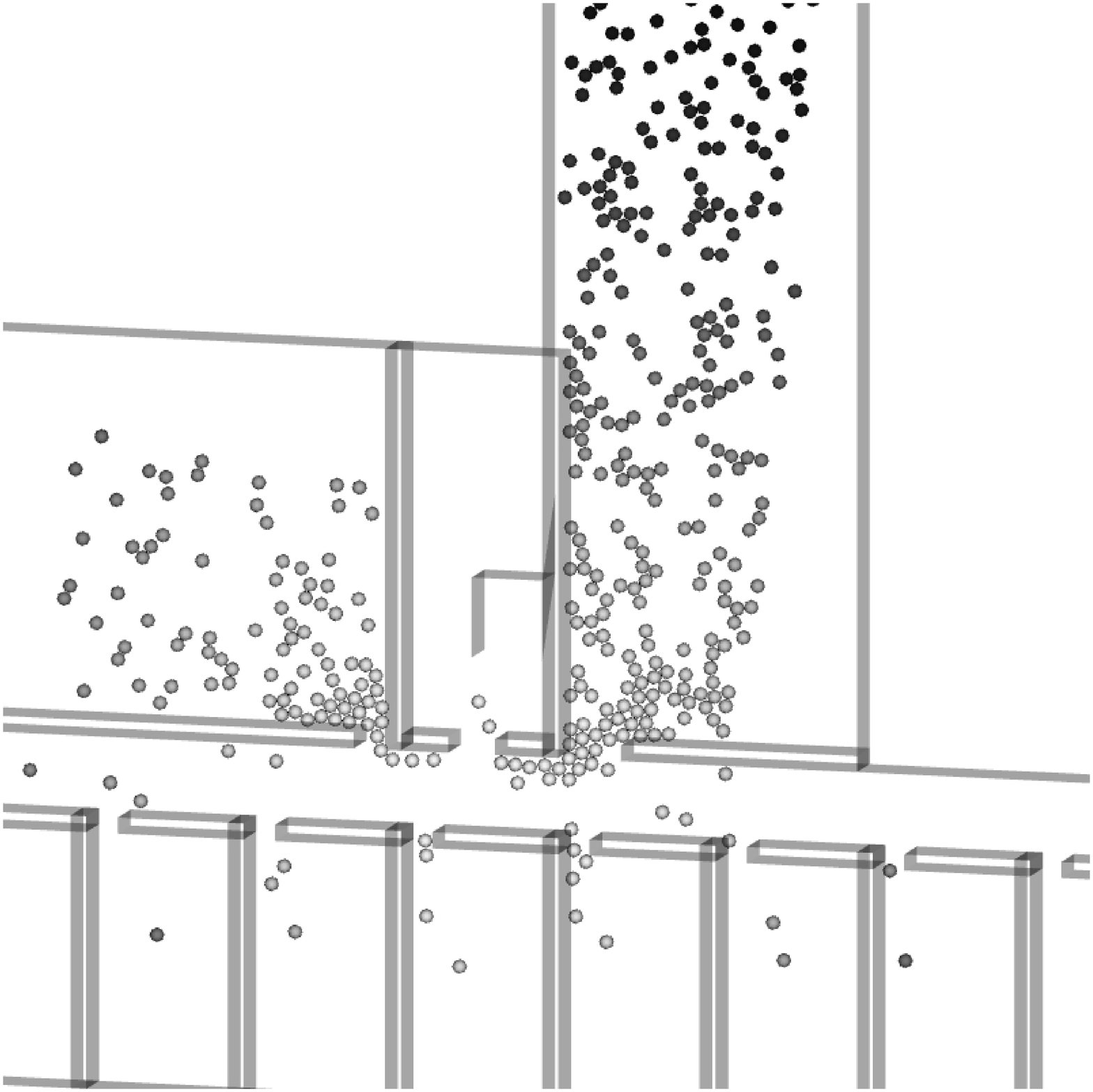}}
\end{tabular}
\vspace{0.8cm}

\begin{tabular}{c@{\hspace{1cm}}c}
\resizebox{!}{!}{\includegraphics[width=0.4\textwidth]{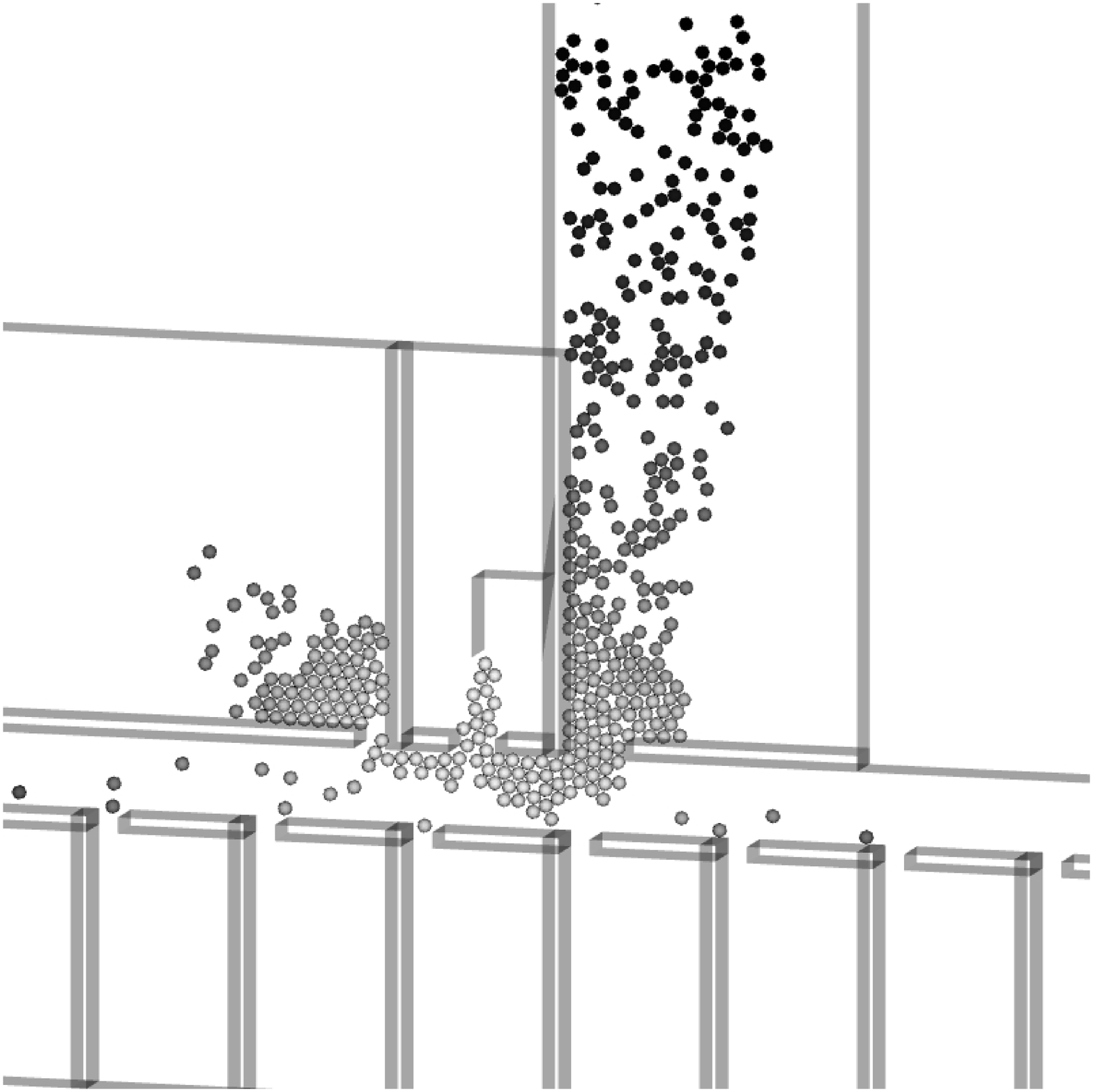}}
&
\resizebox{!}{!}{\includegraphics[width=0.4\textwidth]{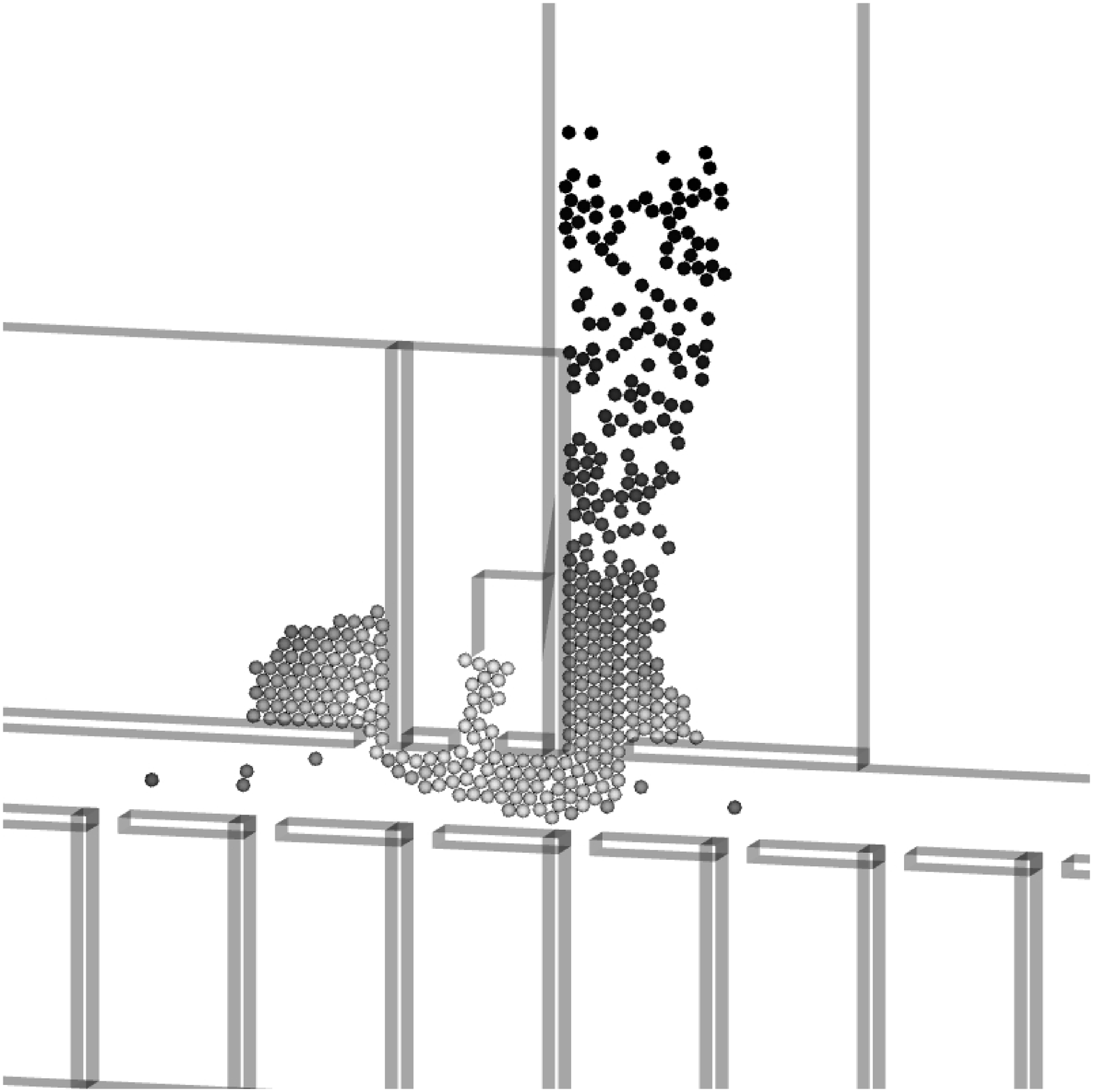}}
\end{tabular}
\vspace{0.8cm}

\begin{tabular}{c@{\hspace{1cm}}c}
\resizebox{!}{!}{\includegraphics[width=0.4\textwidth]{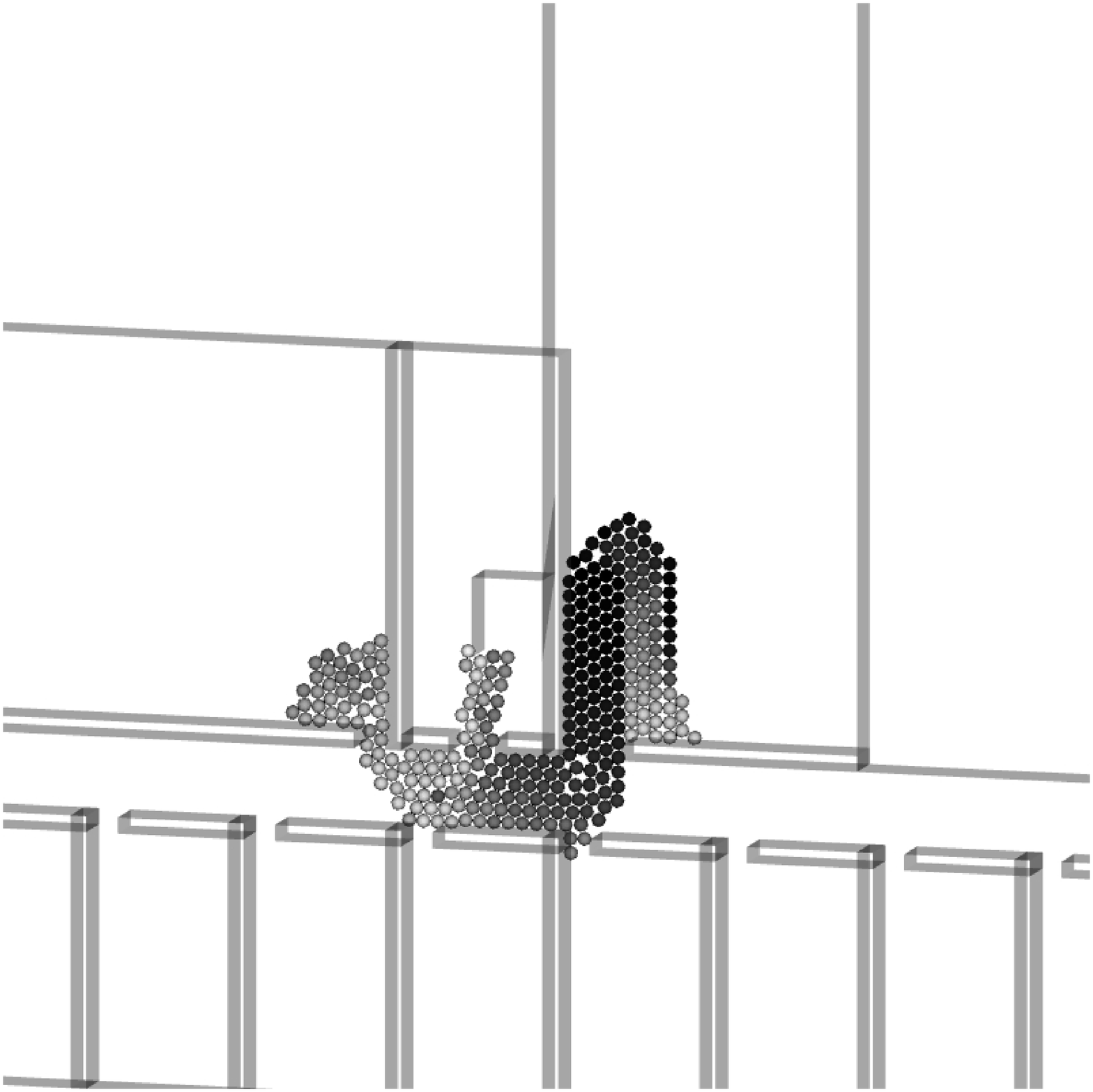}}
&
\resizebox{!}{!}{\includegraphics[width=0.4\textwidth]{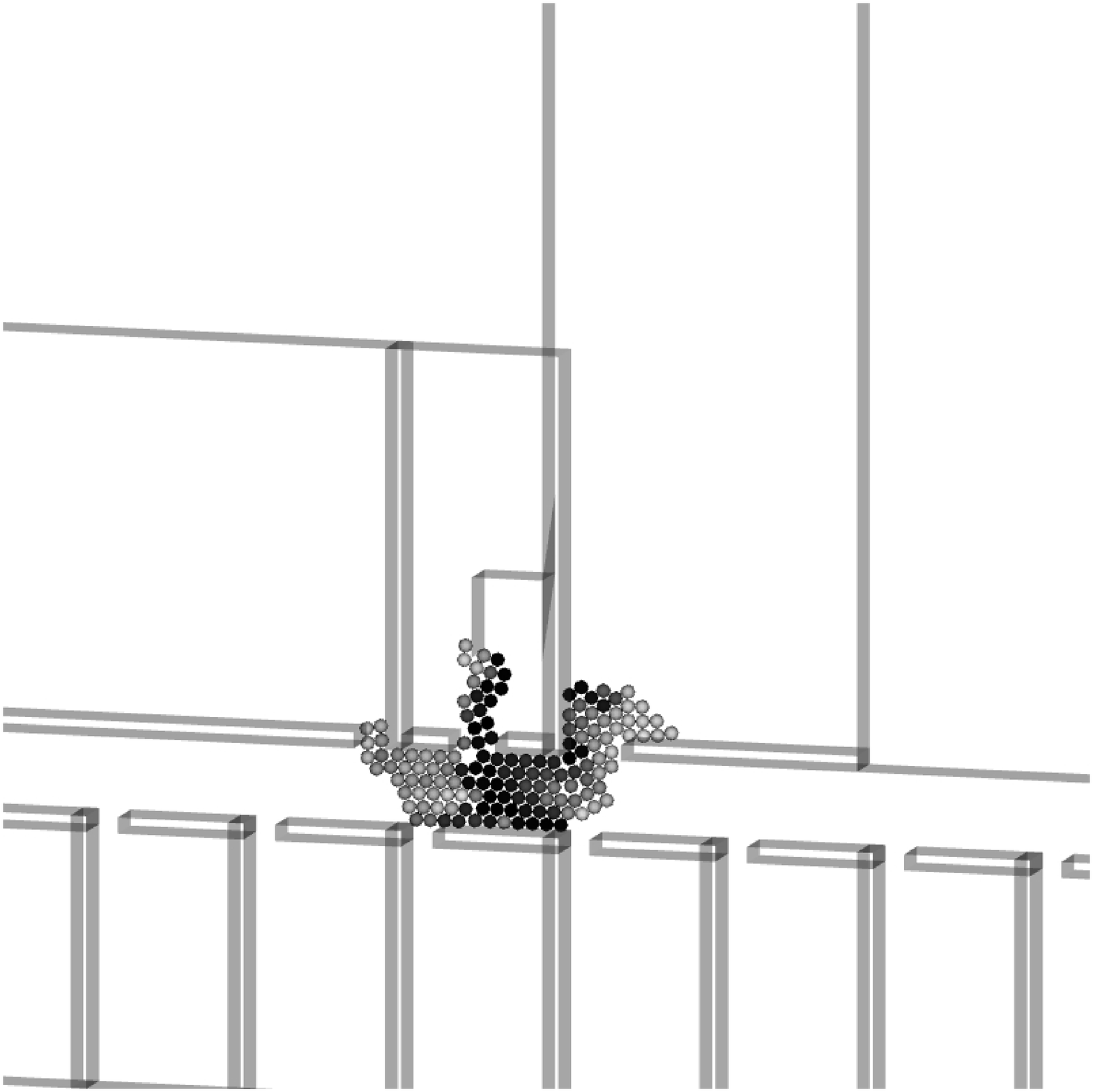}}

\end{tabular}
\end{center}
\caption{Zoom.}
\label{fig:evacbat}
\end{figure}

In the second example  we consider the evacuation of a floor through exit stairs. A zoom on the geometry near the exit (together with the isovalues of the geodesic distance function, on which the spontaneous velocity is built) is represented on Figure~\ref{fig:batmath}. Figure~\ref{fig:evacbat} corresponds to snapshots at times 0s, 5s, 11s, 16s, 41s and 75s.
 Disks are colored according to their initial geodesic distance to the exit. Note that initial ordering is not preserved during the evacuation. Notice also how a jam forms between snapshots 2 and 3 in the room located on the left hand side. This jam decreases significantly the rate at which people exit the room, but it disappears eventually. The final evacuation time is 109s, to be compared to 48s which corresponds to the evacuation time without congestion. 

\bibliographystyle{plain}
\bibliography{biblio}
\end{document}